\documentstyle[amscd,amssymb,verbatim]{amsart}
\newcommand{\lrar}[1]{\begin{picture}(50,10)(-25,-5)                          
\put(-25,0){\vector(1,0){50}}
\put(0,5){\makebox(0,0)[b]{\mbox{$#1$}}}
\end{picture}}

\newcommand{\ldar}[1]{\begin{picture}(10,50)(-5,-25)
\put(0,25){\vector(0,-1){50}}
\put(5,0){\mbox{$#1$}}
\end{picture}}

\newcommand{\ev}{\operatorname{ev}}
\newcommand{\ad}{\operatorname{ad}}

\newcommand{\pr}{\operatorname{pr}}

\newcommand{\coker}{\operatorname{coker}}
\newcommand{\DD}{{\cal D}}

\newcommand{\MM}{{\cal M}}

\newcommand{\SL}{\operatorname{SL}}
\newcommand{\splin}{\operatorname{sl}}
\newcommand{\pgl}{\operatorname{pgl}}
\newcommand{\gl}{\operatorname{gl}}
\newcommand{\G}{{\Bbb G}}

\newcommand{\gog}{{\frak g}}

\newcommand{\lan}{\langle}
\newcommand{\ran}{\rangle}

\newcommand{\CC}{{\cal C}}

\newcommand{\Mat}{\operatorname{Mat}}

\newcommand{\tr}{\operatorname{tr}}

\renewcommand{\P}{{\Bbb P}}

\newcommand{\Pic}{\operatorname{Pic}}

\newcommand{\Ad}{\operatorname{Ad}}
\newcommand{\ga}{\gamma}

\newcommand{\eps}{\epsilon}

\numberwithin{equation}{section}

\newtheorem{thm}{Theorem}
\newtheorem{prop}{Proposition}[section]
\newtheorem{lem}[prop]{Lemma}

\newenvironment{rems}{\vspace{3mm}
\noindent {\bf Remarks.}}{\vspace{3mm}}

\newcommand{\Pf}{\noindent {\it Proof}}
\newcommand{\id}{\operatorname{id}}

\newcommand{\ov}{\overline}

\renewcommand{\Im}{\operatorname{Im}}
\newcommand{\Aut}{\operatorname{Aut}}
\newcommand{\rk}{\operatorname{rk}}
\newcommand{\ra}{\rightarrow}

\newcommand{\PP}{{\cal P}}

\newcommand{\LL}{{\cal L}}
\newcommand{\OO}{{\cal O}}

\newcommand{\Hom}{\operatorname{Hom}}
\newcommand{\Ext}{\operatorname{Ext}}
\newcommand{\End}{\operatorname{End}}
\newcommand{\Res}{\operatorname{Res}}

\renewcommand{\a}{\alpha}
\renewcommand{\b}{\beta}
\newcommand{\om}{\omega}

\newcommand{\la}{\lambda}
\newcommand{\th}{\theta}
\newcommand{\C}{{\Bbb C}}

\newcommand{\Z}{{\Bbb Z}}

\newcommand{\wt}{\widetilde}

\newcommand{\sign}{\operatorname{sign}}

\newcommand{\ed}{\qed\vspace{3mm}}

\title{Classical Yang-Baxter equation and the $A_{\infty}$-constraint}
\author{A. Polishchuk}
\thanks{This work was partially supported by NSF grant}

\begin{document}

\begin{abstract} We show that elliptic solutions of classical
Yang-Baxter equation (CYBE) can be obtained
from triple Massey products on elliptic curve. We  
introduce the associative version of this equation 
which has two spectral parameters and construct its elliptic solutions.
We also study some degenerations of these solutions.
\end{abstract}

\maketitle

\bigskip

\centerline{\sc Introduction}

\bigskip

Recall that the classical Yang-Baxter equation (CYBE) is the 
equation
$$[r^{12}(x),r^{23}(y)]+
[r^{12}(x),r^{13}(x+y)]+
[r^{13}(x+y),r^{23}(y)]=0
$$
where $r(x)$ is a meromorphic function of one complex
variable $x$ in the neighborhood of $0$ taking values in $\gog\otimes
\gog$ for some Lie algebra $\gog$. Here $r^{12}(x)$ denotes the element
$r(x)\otimes 1\in U(\gog)\otimes U(\gog)\otimes U(\gog)$, etc.
In their remarkable paper \cite{BD} Belavin and Drinfeld studied
non-degenerate solutions of the CYBE (i.e. solutions such that
the tensor $r(x)$ has maximal rank for generic $x$) for a simple
Lie algebra $\gog$. They proved that any such solution is equivalent to
either
{\it elliptic, trigonometric, or rational} meaning the character
of dependence of $r(x)$ on $x$. Furthermore, they completely
classified elliptic solutions (which can appear
only in the case $\gog=\splin_n$) and trigonometric solutions. 

In this paper we present an unexpected connection between the CYBE and
the $A_{\infty}$-constraint. The latter is certain generalization of
the associativity axiom invented by Stasheff \cite{S}.
One can consider the notion of 
$A_{\infty}$-algebra (resp. $A_{\infty}$-category)
as a natural replacement for the notion of associative algebra
(resp. category) in the presence of a differential. One of the reasons for
introducing this notion is that the category of dg-algebras
(in which the usual associativity constraint is imposed) doesn't
have enough morphisms, so it is often convenient to embed it into
the larger category of $A_{\infty}$-algebras. 
In this paper we observe that in some special situations triple products
in $A_{\infty}$-category 
\footnote{One technical detail concerning the above relation between
$A_{\infty}$-constraint and the CYBE is that we need to consider
$A_{\infty}$-structures which have {\it cyclic symmetry}.
This notion is defined in \cite{P-hmc} and in \cite{P-ainf} we
showed that there is a cyclic symmetry on the $A_{\infty}$-category
associated with a complex compact manifold.}
can be arranged into tensors satisfying
CYBE. More precisely, we show that all non-degenerate
elliptic solutions of the CYBE for $\splin_n$ arise in this way from certain
triple products in the $A_{\infty}$-version of the derived category
of coherent sheaves on elliptic curve. We also show that all
non-degenerate trigonometric solutions of the CYBE for $\splin_2$
arise in the same way from the $A_{\infty}$-category associated with
the union of two $\P^1$'s glued in two points. 
We expect that one can obtain all non-degenerate
trigonometric solutions of the CYBE for $\splin_n$ by considering
$A_{\infty}$-categories of singular curves of arithmetic genus $1$.

The triple products in $A_{\infty}$-categories leading to CYBE
appear to be specializations of triple products of a more general kind
which in turn produce solutions of another equation that we call
the {\it associative Yang-Baxter equation} (AYBE):
\begin{equation}\label{AYBE}
r^{12}(-u',v)r^{13}(u+u',v+v')-
r^{23}(u+u',v')r^{12}(u,v)+
r^{13}(u,v+v')r^{23}(u',v')=0,
\end{equation}
where $r(u,v)$ is a meromorphic function of two complex variables $(u,v)$
in the neighborhood of $(0,0)$ taking values in $A\otimes A$ where
$A$ is an associative algebra with unit. 
\footnote{Constant solutions of this
equation were considered in \cite{A}.} We conjecture that for $A=\Mat(n,\C)$
the analogue of Belavin-Drinfeld classification holds, i.e.
all non-degenerate solutions of the AYBE are equivalent to either 
elliptic\footnote{Unlike the case of CYBE, ``elliptic'' here means
``elliptic of the third kind'', i.e. we allow functions corresponding
to meromorphic sections of line bundles on an elliptic curve.},
or trigonometric, or rational solutions. In section \ref{scalar} 
we check that this is true
for scalar solutions, i.e. for $A=\C$. 
In this case the only solution is the Kronecker's
function $F(u,v,\tau)$ (see section \ref{KWsec}) and its degenerations.
The relation between AYBE and CYBE is the following.
Let $\pr:\Mat(n,\C)\ra\splin_n(\C)$ be the projection along scalar
matrices. It turns out that in the situations we consider the
function $(\pr\otimes\pr)(r(u,v))$ has a limit as $u\ra 0$. We show that
if $r(u,v)$ satisfies the AYBE and the unitarity condition
\begin{equation}\label{skewsymcond}
r^{21}(-u,-v)=-r(u,v)
\end{equation}
then the limit $\ov{r}(v)=(\pr\otimes\pr)(r(u,v))|_{u=0}$
is a solution of the CYBE. We construct elliptic solutions of the AYBE
for $\Mat(n,\C)$ which specialize in this way
to the usual elliptic $r$-matrices.
Also we construct two trigonometric solutions of the AYBE
for $\Mat(2,\C)$ which
specialize to two different trigonometric solutions of the CYBE
for $\splin_2(\C)$.
In section \ref{reconstr-sec} 
we show that if $\ov{r}(v)$ is a non-degenerate unitary solution of the CYBE
with values in $\splin_n(\C)$ which has no infinitesimal symmetries then 
up to rescaling $r(u,v)\mapsto\exp(cuv)r(u,v)$ (where $c\in\C$)
there exists at most
one unitary solution of the AYBE with values in $\Mat_n(\C)$
of the form $r(u,v)=\frac{1\otimes 1}{u}+r_0(v)+\ldots$ with
$(\pr\otimes\pr)(r_0(v))=\ov{r}(v)$.
This applies in particular to elliptic $r$-matrices since they have
no infinitesimal symmetries.

\noindent {\it Acknowledgment.} I am grateful to Pavel Etingof for
useful discussions, especially for help with proofs of 
Theorems \ref{scthm} and \ref{reconstr}.


\section{Identities between triple Massey products and $r$-matrices}
\label{Massey}

\subsection{Massey products in $A_{\infty}$-categories and
in triangulated categories}
\label{gennon}

Recall that an $A_{\infty}$ category consists of a class of objects, a
collection of (graded) vector spaces of morphisms between them 
equipped with operations $m_n(a_1,\ldots,a_n)$ which associate to any 
sequence $a_1,\ldots, a_n$ of composable morphisms ($n\ge 1$) a new morphism
(of degree $\sum_i \deg(a_i)+2-n$). These operations should satisfy the
set of equations similar to the associativity equations which we call
$A_{\infty}$-{\it constraint}. They have form
$$\sum \pm m_k(a_1,\ldots, a_i, m_l(a_{i+1},\ldots, a_{i+l}),\ldots,a_n)=0$$
where $a_1,\ldots,a_n$ is a sequence of composable morphisms,
the sum is taken over all subsegments in the segment of integers $[1,n]$.
The choice of signs is rather subtle (and non-unique). We follow the
sign convention of \cite{GJ}. For more details regarding this definition
see \cite{P-hmc}. We always impose the condition that our
$A_{\infty}$-category has
strict identity morphisms, i.e. $m_1$-closed elements 
$\id_X\in\Hom^0(X,X)$ for every object $X$, which are units with respect
to $m_2$ and such that any higher product $m_n$ ($n\ge 3$) which has 
$\id_X$ as one of the arguments vanishes.

Loosely speaking Massey products in $A_{\infty}$-categories are expressions
in $m_n$'s which are invariant under arbitrary homotopy of 
$A_{\infty}$-structure (see \cite{P-hmc} for the definition).
Unfortunately, the corresponding formalism seems to be absent in
the existing literature except in the particular case of a
differential graded category which can be considered as an 
$A_{\infty}$-category with $m_n=0$ for $n>2$.
  
On the other hand, there is a definition of Massey products in 
triangulated categories (see \cite{GM} IV.2, \cite{P-Mas}).
These products coincide with the differential graded Massey products
in the case when the triangulated category $\DD$ is {\it enhanced} in the sense
of Bondal-Kapranov's paper \cite{BK}. By definition this means that
$\DD$ is obtained by taking cohomology of a {\it pretriangulated dg-category} 
(the property of a dg-category to be pretriangulated means that
certain convolutions exist). Note that according to Kontsevich's philosophy
(see \cite{K}, \cite{K-t})
this pretriangulated dg-category should be considered as a primary object
(considered up to $A_{\infty}$-equivalence).

The enhanced
triangulated category we are interested in is $D^b(X)$ --- the 
bounded derived
category of coherent sheaves on a projective variety $X$ over a field
$k$ (see \cite{BK}). Let us denote by $D^b_{dg}(X)$ the corresponding 
pretriangulated dg-category. The objects of $D^b_{dg}(X)$ are bounded
complexes of coherent sheaves while the morphisms are given by some standard
complexes computing the corresponding $\Ext$'s. According to general
principles of homological perturbation theory (see \cite{Kad}, \cite{GS},
\cite{GLS},\cite{Markl}) there exists an $A_{\infty}$-category 
$D^b_{\infty}(X)$ with the same objects as $D^b_{dg}(X)$ such that
$D^b_{\infty}(X)$ is $A_{\infty}$-equivalent to $D^b_{dg}(X)$ and
$m_1=0$ in $D^b_{\infty}(X)$. Then Massey products in $D^b(X)$ (as in
triangulated category) and in $D^b_{\infty}$ (as in $A_{\infty}$-category) 
are the same. The advantage of considering $D^b_{\infty}$ is that
we can apply $A_{\infty}$-constraint to derive some non-trivial
relations between Massey products. On the other hand, Massey products in
triangulated categories are easier to compute and they often have
a geometric interpretation. 

In this paper we will only consider triple Massey products of the particular
kind. First, let us recall the definition in the context of triangulated 
categories. Let $X,Y,Z,T$ be objects of a triangulated category $\DD$,
$f\in\Hom(X,Y)$, $g\in\Hom^1(Y,Z):=\Hom(Y,Z[1])$, $h\in\Hom(Z,T)$ be
morphisms, such that $g\circ f=0$, $h\circ g=0$. Then the Massey product
$$MP(f,g,h)\in\coker(\Hom(X,Z)\oplus\Hom(Y,T)
\stackrel{(h,f)}{\ra}\Hom(X,T))$$
is defined as follows. Let
$$Z\stackrel{\a}{\ra} C\stackrel{\b}{\ra} Y\stackrel{g}{\ra} Z[1]\ra\ldots$$
be a distinguished triangle.
Then by assumption there exist morphisms $\wt{f}\in\Hom(X,C)$ and
$\wt{h}\in\Hom(C,T)$ such that
$$\b\circ\wt{f}=f,$$
$$\wt{h}\circ\a=h.$$
The Massey product $MP(f,g,h)$ is defined as the class of the element
$$\wt{h}\circ\wt{f}\in\Hom(X,T).$$

Now let us give a definition of the corresponding triple Massey products in
the context of $A_{\infty}$-categories (see \cite{Fuk}).
Let $X,Y,Z,T$ be objects in an $A_{\infty}$-category $\CC$.
Let us denote by $H\CC$ the graded category obtained from $\CC$
by taking cohomologies of $\Hom$ with  respect to $m_1$.
Then for every triple of morphisms $f\in\Hom^i_{H\CC}(X,Y)$,
$g\in\Hom^j_{H\CC}(Y,Z)$, $h\in\Hom^k_{H\CC}(Z,T)$ such that
$g\circ f=0$, $h\circ g=0$ we can define their Massey
product
$$MP(f,g,h)\in\coker(\Hom^{i+j-1}_{H\CC}(X,Z)\oplus\Hom^{j+k-1}_{H\CC}(Y,T)
\stackrel{(h,f)}{\ra}\Hom^{i+j+k-1}_{H\CC}(X,T)).$$
For this we choose $m_1$-closed elements $\wt{f}\in\Hom^i_{\CC}(X,Y)$,
$\wt{g}\in\Hom^j_{\CC}(Y,Z)$, $\wt{h}\in\Hom^k_{\CC}(Z,T)$
representing $f$, $g$ and $h$. Furthermore, by assumption we have
$$m_2(\wt{f},\wt{g})=m_1(p),$$
$$m_2(\wt{g},\wt{h})=m_1(q)$$
for some $p\in\Hom^{i+j-1}_{\CC}(X,Z)$, $q\in\Hom^{j+k-1}_{\CC}(Y,T)$.
Then we define $MP(f,g,h)$ as the class of the $m_1$-closed element
$$m_3(\wt{f},\wt{g},\wt{h})-m_2(p,\wt{h})+(-1)^{\deg f}m_2(\wt{f},q)$$
(the fact that it is $m_1$-closed follows from the $A_{\infty}$-constraint).
When $m_3=0$ this definition coincides with the usual definition
given in dg-context. On the other hand, if $m_1=0$ then this
Massey product coincides with $m_3$. Finally, we claim 
that this Massey product is
preserved under any equivalence of $A_{\infty}$-categories.
This is a consequence of the following result.

\begin{prop}\label{hominv} 
Let $F:\CC\ra\CC'$ be an $A_{\infty}$-functor
between $A_{\infty}$-categories, let $HF:H\CC\ra H\CC'$
be the induced functor between the corresponding graded categories.
Then
$$HF(MP(f,g,h)=MP(HF(f),HF(g),HF(h)).$$
\end{prop}

\Pf . Let $F=(F_n)$ where $F_n$ are maps from $n$-tuples of
composable morphisms in $\CC$ to morphisms in $\CC'$.
According to the definition of the $A_{\infty}$-functor we have
$$m_2(F_1\wt{f},F_1\wt{g})=F_1m_2(\wt{f},\wt{g})-m_1F_2(\wt{f},\wt{g})=
F_1m_1(p)-m_1F_2(\wt{f},\wt{g}).$$
Since $F_1$ commutes with $m_1$ we get
$$m_2(F_1\wt{f},F_1\wt{g})=m_1(F_1(p)-F_2(\wt{f},\wt{g})).$$
Similarly,
$$m_2(F_1\wt{g},F_1\wt{h})=m_1(F_1(q)-F_2(\wt{g},\wt{h})).$$
Thus, the triple Massey product $MP(HF(f),HF(g),HF(h))$ is
represented by the element
$$
m_3(F_1\wt{f},F_1\wt{g},F_1\wt{h})-m_2(F_1(p)-F_2(\wt{f},\wt{g}),F_1\wt{h})
+(-1)^{\deg f}m_2(F_1\wt{f},F_1(q)-F_2(\wt{g},\wt{h})).$$
Using the identity
\begin{align*}
&
m_3(F_1\wt{f},F_1\wt{g},F_1\wt{h})+m_2(F_2(\wt{f},\wt{g}),F_1\wt{h})
-(-1)^{\deg f}m_2(F_1\wt{f},F_2(\wt{g},\wt{h}))=\\
&F_1m_3(\wt{f},\wt{g},\wt{h})-F_2(m_2(\wt{f},\wt{g}),\wt{h})-
(-1)^{\deg f} F_2(\wt{f},m_2(\wt{g},\wt{h}))
-m_1F_3(\wt{f},\wt{g},\wt{h}).
\end{align*}
we can rewrite the element representing $MP(HF(f),HF(g),HF(h))$ as follows:
\begin{equation}\label{Mas-el}
\begin{array}{l}
F_1m_3(\wt{f},\wt{g},\wt{h})-m_2(F_1(p),F_1\wt{h})+
(-1)^{\deg f}m_2(F_1\wt{f},F_1(q))\\
-F_2(m_1(p),\wt{h})-(-1)^{\deg f} F_2(\wt{f},m_1(q))
-m_1F_3(\wt{f},\wt{g},\wt{h}).
\end{array}
\end{equation}
Note that the last term is a coboundary, hence, it can be omitted.
On the other hand, we have
$$m_2(F_1p,F_1\wt{h})\equiv F_1m_2(p,\wt{h})-F_2(m_1(p),\wt{h})\mod \Im(m_1)$$
and 
$$m_2(F_1\wt{f},F_1q)\equiv F_1m_2(\wt{f},q)+F_2(\wt{f},m_1(q))\mod\Im(m_1).$$
Substituting this in (\ref{Mas-el}) we obtain
that $MP(HF(f),HF(g),HF(h))$ is represented by
$$F_1m_3(\wt{f},\wt{g},\wt{h})-F_1m_2(p,\wt{h})+
(-1)^{\deg f}F_1m_2(\wt{f},q).$$
Therefore, it coincides with $HF(MP(f,g,h))$.
\ed

Both the definitions above can be slightly generalized:
instead of considering a decomposable tensor $f\otimes g\otimes h$
one can take any tensor in the appropriate subspace of
$\Hom^i(X,Y)\otimes\Hom^j(Y,Z)\otimes\Hom^k(Z,T)$.
We leave this to the reader (in the context of triangulated categories
the corresponding definition can be found in \cite{P-Mas}).

\subsection{Generic identity and the associative Yang-Baxter equation}
\label{AYBsec}

Let $\CC$ be an $A_{\infty}$-category with $m_1=0$. 
Assume that we have two
families $\MM$ and $\MM'$ of objects of $\CC$ with the following
properties:

\noindent
(i) for every pair of distinct objects $X_1,X_2\in \MM$
(resp. $Y_1,Y_2\in \MM'$) one has $\Hom^{\bullet}(X_1,X_2)=0$
(resp. $\Hom^{\bullet}(Y_1,Y_2)=0$);

\noindent
(ii) for every $X\in\MM$ and every $Y\in\MM'$ the space
$\Hom^{\bullet}(X,Y)$ is concentrated in degree $0$, the
space $\Hom^{\bullet}(Y,X)$ is concentrated in degree $1$ and 
a perfect pairing
$$\lan\cdot,\cdot\ran:\Hom^0(X,Y)\otimes\Hom^1(Y,X)\ra k$$
is given.

In this situation we can consider the triple products
$$m_3:\Hom^0(X_1,Y_1)\otimes\Hom^1(Y_1,X_2)\otimes\Hom^0(X_2,Y_2)
\ra\Hom^0(X_1,Y_2)$$
and
$$m_3:\Hom^1(Y_1,X_2)\otimes\Hom^0(X_2,Y_2)\otimes\Hom^1(Y_2,X_1)
\ra\Hom^1(Y_1,X_1)$$
where $X_1,X_2\in\MM$, $X_1\neq X_2$, $Y_1,Y_2\in\MM'$, $Y_1\neq Y_2$.
Using the vanishing of the spaces
$\Hom^{\bullet}(X_1,X_2)$ and $\Hom^{\bullet}(Y_1,Y_2)$ and the condition
$m_1=0$ one can immediately see that the corresponding
Massey products coincide with $m_3$.
We assume in addition that the pairing from (ii) is compatible with
these triple products in the following sense:

\noindent
(iii) for every $f_1\in\Hom^0(X_1,Y_1)$, $g_1\in\Hom^1(Y_1,X_2)$,
$f_2\in\Hom^0(X_2,Y_2)$, $g_2\in\Hom^1(Y_2,X_1)$ one has
$$\lan m_3(f_1,g_1,f_2),g_2\ran=-\lan f_1,m_3(g_1,f_2,g_2)\ran=
-\lan m_3(f_2,g_2,f_1),g_1\ran.$$

Note that the condition (iii) is satisfied when $\CC$ has a structure of
{\it cyclic} $A_{\infty}$-category in the sense of \cite{P-ainf}.

Using the duality from (ii) we can rewrite the tensor corresponding to
$m_3$ as a linear map
$$r^{X_1X_2}_{Y_1Y_2}:\Hom^0(X_1,Y_1)\otimes\Hom^0(X_2,Y_2)\ra\Hom^0(X_2,Y_1)
\otimes\Hom^0(X_1,Y_2).$$

\begin{thm} For any triples of distinct objects $X_1,X_2,X_3\in\MM$,
$Y_1,Y_2,Y_3\in\MM'$ one has
\begin{equation}\label{Rid}
(r^{X_3X_2}_{Y_1Y_2})^{12}(r^{X_1X_3}_{Y_1Y_3})^{13}-
(r^{X_1X_3}_{Y_2Y_3})^{23}(r^{X_1X_2}_{Y_1Y_2})^{12}+
(r^{X_1X_2}_{Y_1Y_3})^{13}(r^{X_2X_3}_{Y_2Y_3})^{23}=0
\end{equation}
as a map 
$$\Hom^0(X_1,Y_1)\Hom^0(X_2,Y_2)\Hom^0(X_3,Y_3)\ra
\Hom^0(X_2,Y_1)\Hom^0(X_3,Y_2)\Hom^0(X_1,Y_3).$$
In addition the following skew-symmetry holds:
\begin{equation}
(r^{X_1X_2}_{Y_1Y_2})^{21}=-r^{X_2X_1}_{Y_2Y_1}.
\end{equation}
\end{thm}

\Pf . The skew-symmetry follows easily from the property (iii).
Using it we can rewrite the equation (\ref{Rid}) as follows
$$(r^{X_1X_3}_{Y_2Y_3})^{23}(r^{X_1X_2}_{Y_1Y_2})^{12}+c.p.=0$$
where ``c.p.'' stands for the terms obtained from the first one by cyclic
permutation of indices.
    
Let us consider any six elements 
$f_i\in\Hom^0(X_i,Y_i)$, $g_i\in\Hom^1(Y_i,X_{i+1})$,
where $i\in\Z/3\Z$ (so that $X_4:=X_1$). 
The definition of $r^{X_1X_2}_{Y_1Y_2}$ is equivalent to the following
formula:
$$\lan r^{X_1X_2}_{Y_1Y_2}(f_1\otimes f_2),g_1\ran_1=m_3(f_1,g_1,f_2),$$
where $\lan ?,?\ran_1$ denotes the result of applying the pairing
$\lan ?,?\ran$ in the first component of the tensor product.
It follows that
$$\lan (r^{X_1X_3}_{Y_2Y_3})^{23}(r^{X_1X_2}_{Y_1Y_2})^{12}
(f_1\otimes f_2\otimes f_3),g_1\otimes g_2\ran_{12}=
m_3(m_3(f_1,g_1,f_2),g_2,f_3)$$
where $\lan ?,?\ran_{12}$ denotes the pairing $\lan ?,?\ran$
applied in the first two
components of the tensor product.
Thus, we have
$$\lan (r^{X_1X_3}_{Y_2Y_3})^{23}(r^{X_1X_2}_{Y_1Y_2})^{12}
(f_1\otimes f_2\otimes f_3),g_1\otimes g_2\otimes g_3\ran=
\lan m_3(m_3(f_1,g_1,f_2),g_2,f_3),g_3\ran.$$
Using property (iii) we can rewrite this formula as follows:
\begin{equation}
\lan (r^{X_1X_3}_{Y_2Y_3})^{23}(r^{X_1X_2}_{Y_1Y_2})^{12}
(f_1\otimes f_2\otimes f_3),g_1\otimes g_2\otimes g_3\ran=
-\lan m_3(f_1,g_1,f_2),m_3(g_2,f_3,g_3)\ran.
\end{equation}

On the other hand, applying the $A_{\infty}$-constraint
to five composable morphisms $f_1,g_1,f_2,g_2,f_3$ and using property (i)
we get
\begin{equation}\label{constraint1}
m_3(m_3(f_1,g_1,f_2),g_2,f_3)+
m_3(f_1,m_3(g_1,f_2,g_2),f_3)-
m_3(f_1,g_1,m_3(f_2,g_2,f_3))=0.
\end{equation}
Pairing this identity with $g_3$ and using property (iii) we get
\begin{equation}\label{constcycl}
\begin{array}{c}
\lan m_3(f_1,g_1,f_2),m_3(g_2,f_3,g_3)\ran+
\lan m_3(f_3,g_3,f_1),m_3(g_1,f_2,g_2)\ran+\\
\lan m_3(f_2,g_2,f_3),m_3(g_3,f_1,g_1)\ran=0.
\end{array}
\end{equation}
\ed

Let $A$ be an associative $k$-algebra with a unit.
For a tensor $r^{X_1X_2}_{Y_1Y_2}\in A\otimes_k A$ depending
on two sets of variables $X_1,X_2\in\MM$, $Y_1,Y_2\in\MM'$
the equation (\ref{Rid}) can be considered as an associative version
of the classical Yang-Baxter equation.
In the case when there is no dependence on variables we obtain the equation
$$r^{12}r^{13}-r^{23}r^{12}+r^{13}r^{23}=0$$
which was considered in \cite{A} in connection with infinitesimal
Hopf algebras.

Now let $k=\C$.
Similar to the case of the usual classical Yang-Baxter equation
it is natural to consider solutions with complex variables $X_i,Y_j$
such that $r=r(u,v)$ is a meromorphic function of $u=X_1-X_2$ and
$v=Y_1-Y_2$ (where $u$ and $v$ vary in the neighborhood of $0$).
Then the equation can be rewritten in the form (\ref{AYBE})
while the skew-symmetry equation becomes the equation (\ref{skewsymcond}).
Using the above theorem we will construct below elliptic solutions
of the AYBE
satisying the condition (\ref{skewsymcond}) with values in the matrix algebra
$\Mat(n,\C)$ which specialize to the standard elliptic $r$-matrices
for $\splin_n(\C)$ as $u$ tends to $0$. This limit procedure works more
generally as follows.
We say that  a solution $r(u,v)$ of the AYBE is {\it unitary} if it satisfies
the equation (\ref{skewsymcond}). Similarly, a 
unitary solution of the CYBE is a solution satisfying
the equation $\ov{r}^{21}(-v)=-\ov{r}(v)$.

\begin{lem}\label{limit}
Let $r(u,v)$ be a unitary solution of the AYBE with values in
$\Mat(n,\C)$.
Let $\pr:\Mat(n,\C)\ra\splin_n(\C)$ be the projection along scalar
matrices. Assume that $(\pr\otimes\pr)(r(u,v))$ has a limit as $u\ra 0$.
Then $\ov{r}(v)=(\pr\otimes\pr)(r(u,v))|_{u=0}$ 
is a unitary solution of the CYBE.
\end{lem}

\Pf . Applying the permutation of the first two factors to the equation
(\ref{AYBE}) and making a change of variables $(v,v')\mapsto
(-v,v+v')$, $(u,u')\mapsto (u',u)$ we obtain
$$r^{21}(-u,-v)r^{23}(u+u',v')-
r^{13}(u+u',v+v')r^{21}(u',-v)+
r^{23}(u',v')r^{13}(u,v+v')=0.
$$
Using the equation (\ref{skewsymcond}) this equation can be rewritten
as follows:
$$-r^{12}(u,v)r^{23}(u+u',v')+
r^{13}(u+u',v+v')r^{12}(-u',v)+
r^{23}(u',v')r^{13}(u,v+v')=0.
$$
Subtracting this equation from (\ref{AYBE}) we get
$$[r^{12}(-u',v),r^{13}(u+u',v+v')]-
[r^{23}(u+u',v'),r^{12}(u,v)]+
[r^{13}(u,v+v'),r^{23}(u',v')]=0.
$$
Finally, applying $\pr\otimes\pr$ and substituting $u=u'=0$ we obtain
that $\ov{r}(v)$ satisfies CYBE.
\ed

There is a natural notion of equivalence for the solutions of
(\ref{Rid}). Namely, if $\varphi^X_Y$ is a function with values in
$A^*$ (invertible elements in $A$) and
$r^{X_1X_2}_{Y_1Y_2}$ is a solution of (\ref{Rid})
then
$$\wt{r}^{X_1X_2}_{Y_1Y_2}=(\varphi^{X_2}_{Y_1}\otimes
\varphi^{X_1}_{Y_2})r^{X_1X_2}_{Y_1Y_2}
(\varphi^{X_1}_{Y_1}\otimes\varphi^{X_2}_{Y_2})^{-1}$$
is also a solution of (\ref{Rid}).
We will call the solutions $\wt{r}$ and $r$ {\it equivalent}.
On the other hand, if $\psi_Y$ is a function with values in $\Aut(A)$
then we can construct a new solution by looking at
$$(\psi_{Y_1}\otimes\psi_{Y_2})r^{X_1X_2}_{Y_1Y_2}.$$
However, in the case of the matrix algebra this doesn't give anything
new since all automorphisms are inner.

It is easy to see that if $r(u,v)$ is a solution of (\ref{AYBE}) then
$$c_1\cdot\exp(c_2uv)\cdot r(u,v)$$
is also a solution for arbitrary constants $c_1\in\C^*$ and $c_2\in\C$.
We will call this operation {\it rescaling} of a solution.

It seems reasonable to conjecture that all unitary
solutions of (\ref{AYBE})
with values in the matrix algebra
satisfying 
the non-degeneracy condition (that the tensor $r(u,v)$
is non-degenerate for generic $u,v$) are equivalent (up to rescaling)
to either
{\it elliptic} or {\it trigonometric} or {\it rational} solution
similar to the Belavin-Drinfeld classification in \cite{BD}.
In section \ref{scalar}
we will check our conjecture in the simplest case $n=1$, i.e. we will 
classify scalar unitary solutions of (\ref{AYBE}). 

\subsection{Classical Yang-Baxter equation}\label{YBsec}

Now we will express the ``limit" of $r^{X_1X_2}_{Y_1Y_2}$ as $X_2$
tends to $X_1$ directly in terms of $A_{\infty}$-structure.  
We will see that in the case $X_1=X_2$ 
the Massey products have smaller domain of definition and smaller range
and that the corresponding tensor satisfies the CYBE.

We still consider an $A_{\infty}$-category $\CC$ with $m_1=0$.
Now assume that we have an object $X$ and a family of objects
$\MM$ in $\CC$, such that the following properties hold:

\noindent
(i)' For every pair of distinct objects $Y_1,Y_2\in \MM$
one has $\Hom^{\bullet}(Y_1,Y_2)=0$;
the spaces $\Hom^0(X,X)$ and $\Hom^1(X,X)$ are one-dimensional,
$\Hom^i(X,X)=0$ for $i\neq 0,1$. 

\noindent
(ii)' for every $Y\in\MM$ the space
$\Hom^{\bullet}(X,Y)$ is concentrated in degree $0$, the
space $\Hom^{\bullet}(Y,X)$ is concentrated in degree $1$ and 
the composition map
$$m_2:\Hom^0(X,Y)\otimes\Hom^1(Y,X)\ra \Hom^1(X,X)\simeq k$$
is a perfect pairing.

In this situation we can consider the Massey product induced by the
triple product
\begin{equation}\label{m3X1X2Y}
m_3:\Hom^0(X,Y_1)\otimes\Hom^1(Y_1,X)\otimes\Hom^0(X,Y_2)\ra\Hom^0(X,Y_2)
\end{equation}
where $Y_1,Y_2\in\MM$, $Y_1\neq Y_2$.
The domain of definition of the corresponding triple Massey product
contains tensors
$\sum_i f_i\otimes g_i\otimes h$ such that
$$\sum_i m_2(f_i,g_i)=0.$$
The value of the Massey product on such a tensor is an element of
$\Hom^0(X,Y_2)$ defined up to addition of a scalar multiple of $h$.
It is more convenient to consider the product (\ref{m3X1X2Y}) as a linear
map
$$\Hom^0(X,Y_1)\otimes\Hom^1(Y_1,X)\ra \End(\Hom^0(X,Y_2)).$$
Then the corresponding Massey product is the map
\begin{equation}\label{indmap}
K_{X,Y_1}\ra\End(\Hom^0(X,Y_2))/k\cdot\id,
\end{equation}
where $K_{X,Y_i}\subset\Hom^0(X,Y_i)\Hom^1(Y_i,X)$ is the kernel
of $m_2$.

For every finite-dimensional vector space $V$ over $k$ let us denote by
$\splin(V)\subset\End(V)$ the subspace of traceless endomorphisms,
and $\pgl(V)=\End(V)/k\cdot\id$. We have a canonical isomorphism
$\splin(V)^*\simeq\pgl(V)$ induced by self-duality of $\End(V)$.

Let us choose a linear isomorphism $\tr:\Hom^1(X,X)\ra k$.
Then using the pairing
$$\lan\cdot,\cdot\ran=\tr\circ m_2:\Hom^0(X,Y_i)\otimes\Hom^1(Y_i,X)\ra k$$
we can identify
$\Hom^1(Y_i,X)$ with the dual space to $\Hom^0(X,Y_i)$.
In view of this duality the triple product (\ref{m3X1X2Y}) can be considered
as a tensor
$$\wt{r}_{Y_1,Y_2}\in\End(\Hom^0(X,Y_1))\otimes\End(\Hom^0(X,Y_2)).$$
On the other hand, $K_{X,Y_1}$ can be identified with the subspace
$\splin(\Hom^0(X,Y_1))\subset\End(\Hom^0(X,Y_1))$.
Thus, we can rewrite the map (\ref{indmap}) as a linear map
$$\splin(\Hom^0(X,Y_1))\ra\pgl(\Hom^0(X,Y_2))$$
or equivalently as a tensor
$$r_{Y_1,Y_2}=r^X_{Y_1,Y_2}\in\pgl(\Hom^0(X,Y_1))\otimes
\pgl(\Hom^0(X,Y_2)).$$
It is easy to see that $r_{Y_1,Y_2}$ is the image of $\wt{r}_{Y_1,Y_2}$
under the natural projection.
By Proposition \ref{hominv} 
the tensor $r_{Y_1,Y_2}$ is invariant under any homotopy
of $A_{\infty}$-structure. 

We assume in addition that 

\noindent (iii)'
for every $f_i\in\Hom^0(X,Y_i)$, $g_i\in\Hom^1(Y_i,X)$, $i=1,2$,
one has
$$\lan m_3(f_1,g_1,f_2),g_2\ran=-\lan f_1,m_3(g_1,f_2,g_2)\ran=
-\lan m_3(f_2,g_2,f_1),g_1\ran.$$

\begin{thm}\label{CYBE-thm}
For every triple of distinct objects $Y_1,Y_2,Y_3\in\MM$
one has
\begin{equation}\label{CYBE}
[r_{Y_1,Y_2}^{12},r_{Y_1,Y_3}^{13}]+
[r_{Y_1,Y_2}^{12},r_{Y_2,Y_3}^{23}]+
[r_{Y_1,Y_3}^{13},r_{Y_2,Y_3}^{23}]=0
\end{equation}
in the Lie algebra
$\pgl(\Hom^0(X,Y_1))\otimes\pgl(\Hom^0(X,Y_2))\otimes\pgl(\Hom^0(X,Y_3))$.
In addition the following skew-symmetry holds:
\begin{equation}
r_{Y_1,Y_2}^{21}=-r_{Y_2,Y_1}.
\end{equation}
\end{thm}

\Pf . Let us consider six elements 
$f_i\in\Hom^0(X,Y_i)$, $g_i\in\Hom^1(Y_i,X)$,
where $i\in\Z/3\Z$, such that $\lan f_i,g_i\ran=0$ for all $i$.
In fact, the argument below should (and can) be applied to 
a slightly more general data: each tensor $f_i\otimes g_i$ should
be replaced by an arbitrary element of $K_{X,Y_i}$. However, we
restrict ourself to the case of decomposable tensors to simplify
notations. By definition we have
$$\lan \wt{r}_{Y_1,Y_2}(f_1\otimes f_2),g_1\otimes g_2\ran=
\lan m_3(f_1,g_1,f_2),g_2\ran.$$
Together with the property (iii)' this immediately implies the skew-symmetry
of $r$. Using it we can rewrite the equation (\ref{CYBE}) in the following
form:
$$[r_{Y_1,Y_2}^{12},r_{Y_2,Y_3}^{23}]+c.p.=0.$$
It is easy to see that
\begin{align*}
&\lan \wt{r}_{Y_2,Y_3}^{23}\wt{r}_{Y_1,Y_2}^{12}(f_1\otimes f_2\otimes f_3),
g_1\otimes g_2\otimes g_3\ran=
\lan m_3(m_3(f_1,g_1,f_2),g_2,f_3),g_3\ran=\\
&-\lan m_3(f_1,g_1,f_2), m_3(g_2,f_3,g_3)\ran.
\end{align*}
The $A_{\infty}$-constraint applied to the morphisms $f_1$, $g_1$, $f_2$,
$g_2$, $f_3$ differs from (\ref{constraint1}) by one additional term:
\begin{align*}
&m_3(m_3(f_1,g_1,f_2),g_2,f_3)+
m_3(f_1,m_3(g_1,f_2,g_2),f_3)-
m_3(f_1,g_1,m_3(f_2,g_2,f_3))-\\
&m_2(m_4(f_1,g_1,f_2,g_2),f_3)=0.
\end{align*}
However, this additional term drops out when we apply pairing with $g_3$
since $m_4(f_1,g_1,f_2,g_2)$ is a multiple of $\id_X$ and $\lan f_3,g_3\ran=0$.
Thus, the equality (\ref{constcycl}) still holds in our situation.
It follows that the tensor 
$$\wt{r}_{Y_2,Y_3}^{23}\wt{r}_{Y_1,Y_2}^{12}+c.p.\in
\End(\Hom^0(X,Y_1)\otimes\Hom^0(X,Y_2)\otimes\Hom^0(X,Y_3))$$
is orthogonal to $\splin(\Hom^0(X,Y_1))\otimes\splin(\Hom^0(X,Y_2))\otimes
\splin(\Hom^0(X,Y_3))$.
Hence, its projection to
$$\pgl(\Hom^0(X,Y_1))\otimes\pgl(\Hom^0(X,Y_2)\otimes\pgl(\Hom^0(X,Y_3))$$
is zero.
Similar statement holds for the tensor $\wt{r}^{12}\wt{r}^{23}+c.p.$
so we are done.
\ed

Assuming in addition that all the spaces $\Hom^0(X,Y)$ for $Y\in\MM$
have the same dimension $n$ (this is true in all examples)
we can choose isomorphisms $\Hom^0(X,Y)\simeq k^n$ and consider
$r_{Y_1,Y_2}$ as an element of $\pgl_n\otimes\pgl_n$. Then
the map $(Y_1,Y_2)\mapsto r_{Y_1,Y_2}$ defined on all pairs
such that $Y_1\not\simeq Y_2$ is a solution of the CYBE
for $\pgl_n$. A different choice of isomorphisms $\Hom^0(X,Y)\simeq k^n$
leads to an equivalent solution.
In the case $k=\C$ one often has a situation when
objects $X_i$ and $Y_j$ are parametrized by complex variables and all the
spaces $\Hom(X_i,Y_j)$ can be identified with $\C^n$ in such a way
that tensors $r^{X_1,X_2}_{Y_1,Y_2}$ (resp. $r^X_{Y_1,Y_2}$) depend only
on differences of complex parameters corresponding to $X_1,X_2$ and
$Y_1,Y_2$. In this case the solutions of the CYBE corresponding to
$r^{X}_{Y_1,Y_2}$ are obtained from the solutions of the
AYBE corresponding to $r^{X_1,X_2}_{Y_1,Y_2}$ by the limit procedure
described in lemma \ref{limit}.

The above proof also shows that the tensor $r_{Y_1,Y_2}\in\pgl_n\otimes\pgl_n$
has the following property in addition to the CYBE: there exists
a lifting $\wt{r}_{Y_1,Y_2}\in\gl_n\otimes\gl_n$ of $r_{Y_1,Y_2}$
such that 
$$\wt{r}_{Y_2,Y_3}^{23}\wt{r}_{Y_1,Y_2}^{12}+c.p.$$
projects to zero in $\pgl_n^{\otimes 3}$.
It would be interesting to study which solutions of the CYBE satisfy
this property.

\subsection{Spherical objects}
\label{spher}

Let $\DD$ be a triangulated category over a field $k$, such that
all spaces $\Hom(X,Y)$ are finite-dimensional. We use the notation
$\Hom^i(X,Y):=\Hom(X,Y[i])$.

Following \cite{ST} we call an object $F\in\DD$ $n$-{\it spherical} if
$\Hom^i(F,F)=0$ for $i\neq 0,n$, $\Hom^0(F,F)\simeq\Hom^n(F,F)\simeq k$,
and for every $X\in\DD$ the composition map
$$\Hom^i(F,X)\Hom^{n-i}(X,F)\ra\Hom^n(F,F)\simeq k$$
is a perfect pairing.

In the case when $\DD$ is enhanced in the sense of \cite{BK}
one can define the autoequivalence $T_F:\DD\ra\DD$ such that
for every object $X\in\DD$ with $\Hom^i(F,X)=0$ for $i\neq 0$
there is an exact triangle
$$\Hom^0(F,X)\otimes F\ra X\ra T_F X\ra\ldots$$
The case when $\DD$ is a subcategory in the bounded derived category of
quasicoherent sheaves on a projective variety was considered in
details by Seidel and Thomas in \cite{ST}.
The general case of an enhanced triangulated
category is similar. It seems that the construction
of the functor $T_F$ can be generalized to the case when
$\DD$ has a structure of triangulated
$A_{\infty}$-category as defined by Kontsevich \cite{K-t}.

It is easy to see that all spherical objects in the derived category
of coherent sheaves on an elliptic curve $E$ are (up to shift) either simple
vector bundles or structure sheaves of points. In particular, we
observe that the group of autoequivalences of $D^b(E)$ acts transitively
on the set of isomorphism classes of spherical objects.
It seems to be an interesting problem to classify spherical objects
in the case when $E$ is replaced by a singular projective curve of
arithmetic genus $1$. It is natural to consider only such curves
for which the structure sheaf $\OO$ coincides with the dualizing sheaf.
In this case $\OO$ and structure sheaves of smooth points are spherical.
The corresponding functor $T_{\OO}$ together with tensorings by line
bundles and automorphisms of the curve generate a large group of
autoequivalences of the derived category. In particular, we obtain
a lot of spherical objects. However, it is not known whether in this
case the group of autoequivalences acts transitively on spherical objects.

\subsection{Non-degeneracy criterion}

From now on we will always work in an enhanced triangulated category
which has a cyclic symmetry considered as an $A_{\infty}$-category.
We also keep the notations of sections \ref{YBsec} and \ref{spher}.
Recall that a tensor $t\in V_1\otimes V_2$ is called
{\it non-degenerate} if it
induces an isomorphism $V_1^{\vee}\ra V_2$.
We define the non-degeneracy condition
for the tensor $r^{X_1X_2}_{Y_1Y_2}$ by considering it
as an element of
$$(\Hom^0(X_1,Y_1)^{\vee}\otimes\Hom^0(X_2,Y_1))\otimes
(\Hom^0(X_2,Y_2)^{\vee}\otimes\Hom^0(X_1,Y_2)).$$

\begin{thm}\label{nondeg} 
Assume that $Y_1$ and $Y_2$ are $1$-spherical.
Then the tensor
$r^{X_1X_2}_{Y_1Y_2}$ (resp. $r^X_{Y_1,Y_2}$) is non-degenerate
if and only if $\Hom^i(T_{Y_2}X_1,T_{Y_1}X_2)=0$ (resp.
$\Hom^i(T_{Y_2}X,T_{Y_1}X)=0$) for $i=1,2$.
\end{thm}

\Pf . Let us first consider the tensor $r^X_{Y_1,Y_2}$.
Using the definition of the Massey product in the context of
triangulated categories (see section \ref{gennon})
we obtain that $r^X_{Y_1,Y_2}$
corresponds to the composition map
\begin{equation}\label{compmap}
\Hom^0(X,T_{Y_1}X)\otimes\Hom^0(T_{Y_1}X,Y_2)\ra\Hom^0(X,Y_2).
\end{equation}
More precisely, the exact triangle
$$X\ra T_{Y_1}X\ra\Hom^1(Y_1,X)\otimes Y_1\ra\ldots$$
induces the exact sequence
$$0\ra\Hom^0(X,X)\ra\Hom^0(X,T_{Y_1}X)\ra K_{X,Y_1}\ra 0$$
and an isomorphism
$$\Hom^0(T_{Y_1}X,Y_2)\wt{\ra}\Hom^0(X,Y_2).$$
Thus, we have a commutative diagram
\begin{equation}
\begin{array}{ccc}
\Hom^0(X,T_{Y_1}X) & \stackrel{\a}{\ra}\Hom^0(T_{Y_1}X,Y_2)^{\vee}\Hom^0(X,Y_2)\simeq
&\Hom^0(X,Y_2)^{\vee}\Hom^0(X,Y_2)\\
\setlength{\unitlength}{0.15mm}
\ldar{} & & \setlength{\unitlength}{0.15mm}\ldar{}\\
\setlength{\unitlength}{0.80mm}
K_{X,Y_1} &\lrar{r^X_{Y_1,Y_2}} & \pgl(\Hom^0(X,Y_2))
\end{array}
\end{equation}
where the map $\a$ is obtained from (\ref{compmap}) by dualization.
By definition the map $\a$ sends the one-dimensional subspace
$\Hom^0(X,X)\subset\Hom^0(X,T_{Y_1}X)$ to the span of the identity
in $\End(\Hom^0(X,Y_2))$.
Thus, the tensor $r^X_{Y_1Y_2}$ is non-degenerate if and only if
$\a$ is an isomorphism.
To this end we observe that $\a$ is obtained by applying the functor
$\Hom^0(X,?)$ to the second arrow of the following exact triangle:
$$T_{Y_2}^{-1}T_{Y_1}X\ra T_{Y_1}X\ra\Hom^0(T_{Y_1}X,Y_2)^{\vee}\otimes Y_2
\ra\ldots$$
If $\Hom^i(X,T_{Y_2}^{-1}T_{Y_1}X)=0$ for $i=0,1$ then clearly,
$\a$ is an isomorphism.
To show that the converse is true we have to check that
$\Hom^{-1}(X,Y_2)=0$ and $\Hom^1(X,T_{Y_1}X)=0$.
The first vanishing holds by the assumption (ii)'.
From the exact triangle defining
$T_{Y_1}X$ we obtain the following long exact sequence:
$$\Hom^0(X,Y_1)\Hom^1(Y_1,X)\ra\Hom^1(X,X)\ra\Hom^1(X,T_{Y_1}X)\ra
\Hom^1(X,Y_1)\Hom^1(Y_1,X)\ra\ldots$$
Now the condition (ii)' implies that the first arrow is surjective
and the last term vanishes, hence, $\Hom^1(X,T_{Y_1}X)=0$.

In the case of the tensor $r^{X_1X_2}_{Y_1Y_2}$ the proof is very similar
(but more simple): one has natural isomorphisms
$$\Hom^0(X_1,T_{Y_1}X_2)\simeq\Hom^0(X_1,Y_1)\otimes\Hom^1(Y_1,X_2),$$
$$\Hom^0(T_{Y_1}X_2,Y_2)\simeq\Hom^0(X_2,Y_2),$$
while the corresponding Massey product is given by a composition
$$\Hom^0(X_1,T_{Y_1}X_2)\otimes\Hom^0(T_{Y_1}X_2,Y_2)\ra\Hom^0(X_1,Y_2)$$
Thus, the non-degeneracy is equivalent to the condition that the
map
$$\Hom^0(X_1,T_{Y_1}X_2)\ra\Hom^0(T_{Y_1}X_2,Y_2)^{\vee}\otimes
\Hom^0(X_1,Y_2)$$
is an isomorphism. Now the proof can be completed similar to the case
of $r^X_{Y_1Y_2}$.
\ed

\subsection{Solutions associated with simple vector bundles}
\label{bundle-sec}

Now let us consider a more specific situation in which
the general categorical setup
described above is realized. Namely as an enhanced triangulated category
we will take the derived category of a projective curve $C$ of
arithmetic genus $1$. The objects $X_i$ will be simple vector bundles
while the objects $Y_i$ will be structure sheaves of smooth points.
For simplicity let us assume that
$C$ is reduced and it is either irreducible or
it is a union of $\P^1$'s intersecting transversally.
Then the dualizing sheaf of $C$ is $\OO_C$
which implies that most of the conditions (i)-(iii)
(resp. (i)'-(iii)') are satisfied automatically.
More precisely, to check them one can use
the following two lemmas (which are easy consequences of
Riemann-Roch theorem and Serre duality 
on the curve $C$).  

\begin{lem} Let $V$ be a vector bundle on $C$. Then
$\chi(C,V)=\deg V$ where $\deg(V)$ is the sum of degrees of
restrictions of $V$ to irreducible components of $C$.
\end{lem}

\begin{lem} Let $X$ be a simple vector bundle on $C$ 
or a structure sheaf of a smooth point on $C$.
Then $\Ext^i(X,X)=0$ for $i\neq 0,1$,
$\Ext^1(X,X)\simeq k$ and the pairing
$$\Hom(X,Y)\otimes\Hom(Y,X[1])\ra\Ext^1(X,X)\simeq k$$
is non-degenerate for any object $Y$ of the bounded derived
category of coherent sheaves on $C$.
\end{lem}

The only remaining
condition to be checked is that all $\Hom^0$ and $\Ext^1$ between
two simple bundles in question vanish. For example, this is true when these 
bundles are of the form $(V,V\otimes\LL)$ where $\LL$ is a line bundle on $C$
which has degree zero and is not annihilated by $\rk V$ in $\Pic(C)$.
The corresponding triple Massey products are computed in the following
theorem.

\begin{thm}\label{bundlethm}
(a) Let $V_1$, $V_2$ be a pair of simple bundles on $C$ such that
$\Hom^0(V_1,V_2)=\Ext^1(V_1,V_2)=0$. Let $y_1,y_2$ be a pair of distinct
smooth points of $C$. Then the tensor
$$r^{V_1,V_2}_{\OO_{y_1},\OO_{y_2}}\in V_{1,y_1}\otimes V_{2,y_1}^{\vee}
\otimes V_{1,y_2}^{\vee}\otimes V_{2,y_2}$$ 
corresponds to the following composition
$$\Hom(V_{1,y_1},V_{2,y_1})\lrar{\Res_{y_1}^{-1}}
\Hom(V_1,V_2(y_1))\lrar{\ev_{y_2}}\Hom(V_{1,y_2},V_{2,y_2})$$
where the map
$$\Res_{y}:\Hom(V_1,V_2(y))\wt{\ra}\Hom(V_{1,y},V_{2,y})$$
is obtained by taking the residue at a smooth point $y$,
the map $\ev_y$ is the evaluation at a point $y$.

\noindent
(b) Let $V$ be a simple bundle on $C$. Then the tensor
$$r^{V}_{\OO_{y_1},\OO_{y_2}}\in \splin(V_{y_1})\otimes\splin(V_{y_2})$$
corresponds to the composition
$$\splin(V_{y_1})\lrar{\Res_{y_1}^{-1}}H^0(C,\ad V(y_1))
\lrar{\ev_{y_2}}\splin(V_{y_2})$$
where $\ad V$ is the bundle of traceless endomorphisms of $V$.

\noindent
(c) If $V_2\not\simeq V_1(y_2-y_1)$ (resp. $V\not\simeq V(y_2-y_1)$)
then the tensor
$r^{V_1,V_2}_{\OO_{y_1},\OO_{y_2}}$ in (a) (resp.
$r^V_{\OO_{y_1},\OO_{y_2}}$ in (b)) is non-degenerate.
\end{thm}

\Pf . (a) Let us choose an isomorphism between the dualizing sheaf on $C$
and $\OO_C$. By Serre duality we have
$$\Ext^1(\OO_{y_1},V_2)\simeq
\Hom(V_2,\OO_{y_1})^*\simeq V_{2,y_1}.$$
Moreover, the universal extension sequence
$$0\ra V_2\ra U\ra \Ext^1(\OO_{y_1},V_2)\otimes \OO_{y_1}\ra 0$$
can be identified with the canonical exact sequence
\begin{equation}\label{seq}
0\ra V_2\ra V_2(y_1)\ra V_2(y_1)|_{y_1}\ra 0
\end{equation}
where the isomorphism $\OO(y_1)|_{y_1}\simeq\OO_{y_1}$ is induced
by the trivialization of the dualizing sheaf on $C$.
Now by definition of the triple Massey products in triangulated categories
we have to consider the composition map
$$\Hom(V_1, V_2(y_1))\otimes\Hom(V_2(y_1),\OO_{y_2})\ra\Hom(V_1,\OO_{y_2})$$
and use the isomorphisms
$$\Hom(V_1,V_2(y_1))\wt{\ra}\Hom(V_1,V_2|_{y_1})$$
$$\Hom(V_2(y_1),\OO_{y_2})\wt{\ra}\Hom(V_2,\OO_{y_2})$$
induced by the sequence (\ref{seq}). By definition the first
of these isomorphisms is given by taking the residue at $y_1$,
so we arrive at the required description of the Massey product. 

\noindent (b) The proof is analogous to (a) and is omitted.

\noindent (c) It is known (see \cite{ST}) that for any smooth point
$y\in C$ the object $\OO_y$ is spherical and the corresponding
functor $T_{\OO_y}$ is given by tensoring with the line bundle
$\OO_C(y)$.
Thus, by theorem \ref{nondeg} the tensor $r^{V_1,V_2}_{\OO_{y_1},\OO_{y_2}}$
is non-degenerate if and only if
$$\Ext^i(V_1(y_2),V_2(y_1))=0$$
for $i=0,1$. Note that the Riemann-Roch theorem for vector bundles on $C$
implies that
$$h^1(C,V_1^{\vee}\otimes V_2(y_1-y_2))=
h^0(C,V_1^{\vee}\otimes V_2(y_1-y_2)).$$
Since $V_1$ and $V_2(y_1-y_2)$ are non-isomorphic simple bundles we have
$\Hom(V_1,V_2(y_1-y_2))=0$, therefore $\Ext^1(V_1,V_2(y_1-y_2))=0$.
The case of the tensor $r^V_{\OO_{y_1},\OO_{y_2}}$ is similar.
\ed

Combining this theorem with theorem \ref{CYBE-thm} we obtain non-degenerate
solutions of the AYBE and of the CYBE associated with simple bundles on a 
projective curve $C$ of arithmetic genus $1$ with trivial dualizing sheaf.
More precisely, we also have to choose a connected component $C_0$ of $C$
in which points $y_i$ vary. If we fix a point $y_0\in C$ and a
uniformization of $C_0\cap C^{reg}$ compatible with the group law
on the set of smooth points $C^{reg}$ of $C$, then we can consider
the tensor $r$ as depending on complex parameters (two parameters in the
case of the AYBE and one parameter in case of the CYBE).
It is known that in the case when $C$ is an elliptic
curve one obtains all non-degenerate elliptic solutions of the CYBE by the
procedure described in Theorem \ref{bundlethm} (b).
In section \ref{trig-sec} we will construct a simple bundle of rank $2$
on the union of two $\P^1$'s intersecting in two points.
Considering points on two different components of this curve we will
obtain two different trigonometric solutions of the CYBE for $\splin_2$.
In each of these cases (elliptic and trigonometric for $\splin_2)$
we also construct solutions of the AYBE specializing
to the solutions of the CYBE.

\section{Elliptic solutions}
\label{ell-sec}

\subsection{Non-degenerate elliptic solutions}

Let $E$ be an elliptic curve over a field $k$, 
$V$ be a simple vector bundle on $E$, i.e. such that
$\Hom(V,V)\simeq k$. Note that $V$ is a $1$-spherical
object in the derived category of coherent sheaves on $E$. 
Assume that $V$ has positive degree $d$. Then we can apply
the construction of the tensor
$r^{X_1X_2}_{Y_1Y_2}$ (resp. $r^X_{Y_1,Y_2}$) from section \ref{AYBsec}
(resp. section \ref{YBsec})
to $X_i$ varying in a family of line bundles of degree zero
(resp. $X=\OO_E$), $Y_j$ varying in a family of bundles
obtained from $V$ be translation. Note that this is essentially
equivalent to the situation of section \ref{bundle-sec}
since applying the Fourier-Mukai transform to structure sheaves
of points one gets line bundles of degree $0$.
Let $e\in E$ be the neutral element. We fix a trivialization of 
$\det V$ (top wedge power of $V$) at $e$. 
For every $x\in E(k)$ let us consider the following line bundle on $E$
trivialized at $e$: 
$$\PP^d_x=t_x^*\det V\otimes (\det V)^{-1}\otimes (\det V)^{-1}|_x$$
where $t_x:E\ra E$ is the translation by $x$.
Note that $\PP^d_x$ depends on $V$ only through its degree $d$
which is reflected in the notation. The map
$x\mapsto \PP^d_x$ is a homomorphism from $E(k)$ to the Picard
group of $E$. Furthermore, if we denote 
$$\lan x,y\ran^d=(\PP^d_x)|_y$$
then $\lan ?,?\ran^d$ is a symmetric biextension of $E\times E$.
We claim that there exists a line bundle $L$ on $E$ such that
for every $x\in E(k)$ there is a canonical isomorphism
$$t_{rx}^*V\simeq L|_x\otimes\PP^d_x\otimes V,$$
where $r$ is the rank of $V$. Indeed, since the isomorphism
class of a simple vector bundle is determined by its determinant
it suffices to check that $t_{rx}^*V$ and $\PP^d_x\otimes V$
have the same determinants which is clear (in fact, using the theorem of
the cube one can show that $L\simeq(\det V)^r$).
Thus, for every $x,y\in E(k)$ we have a sequence of isomorphisms
\begin{align*}
&\Hom(\PP^d_x,t_y^*V)\simeq H^0(E,\PP^d_{-x}\otimes t_y^*V)\simeq
\lan x,y\ran^d\otimes H^0(E,\PP^d_{-x}\otimes V)\simeq \\
&\lan x,y\ran^d\otimes L^{-1}|_{-x}\otimes H^0(E,t_{-rx}^*V)\simeq 
\lan x,y\ran^d\otimes L^{-1}|_{-x}\otimes H^0(E,V).
\end{align*}
Thus, the function
$$(x_1,x_2;y_1,y_2)\mapsto r_V(x_1,x_2;y_1,y_2):=
r^{\PP^d_{x_1},\PP^d_{x_2}}_{t_{y_1}^*V,t_{y_2}^*V}$$
takes values in 
$$\lan x_2-x_1,y_1-y_2\ran^d\otimes\End(H^0(E,V))\otimes\End(H^0(E,V))$$ 
while the function
$$(y_1,y_2)\mapsto r_V(y_1,y_2):=r^{\OO}_{t_{y_1}^*V,t_{y_2}^*V}$$
takes values in $\pgl(H^0(E,V))\otimes\pgl(H^0(E,V))$.
Note that $r_V(x_1,x_2;y_1,y_2)$ is defined only
when $\PP^d_{x_1}\not\simeq \PP^d_{x_2}$ and
$t_{y_1}^*V\not\simeq t_{y_2}^*V$ which happens presicely when
$d(x_1-x_2)\neq 0$ and $d(y_1-y_2)\neq 0$ in $E$.
Similarly, $r^V(y_1,y_2)$ is defined for $d(y_1-y_2)\neq 0$ in $E$.
Also it is easy to see that
$r_V(x_1,x_2;y_1,y_2)$ (resp. $r_V(y_1,y_2)$)
actually depends only on the differences
$x_1-x_2$ and $y_1-y_2$ (resp. on $y_1-y_2$). So we will use the
notation
$$r_V(x;y)=r_V(0,x;0,y),$$
$$r_V(y)=r_V(0,y).$$

Now we will show that the non-degeneracy criterion
of theorem \ref{nondeg} applies to these tensors for generic
values of parameters.

\begin{prop} Assume that $x,y\in E(k)$ are such that
$dx\neq 0$, $dy\neq 0$, $d(dy-x)\neq 0$ (resp.
$y\in E(k)$ is such that $d^2y\neq 0$).Then the tensor
$r_V(x,y)$ (resp. $r_V(y)$) is non-degenerate.
\end{prop}

\Pf . Using the action of a central extension of
$\SL_2(\Z)$ of $D^b(E)$ (see \cite{M}, \cite{P-FM})
we can find an
autoequivalence $S:D^b(E)\ra D^b(E)$ which sends 
a pair of bundles $(V,t_y^*V)$ to the pair of sheaves
$(\OO_{y_1},\OO_{y_2})$ for some points $y_1\neq y_2$. Then
$S(\OO_E)$ and $S(\PP^d_x)$ are simple vector bundles of rank $d$.
Since the twist functors $T_{\OO_{y_i}}$ are just tensorings
by $\OO_E(y_i)$ we have only to check that 
$$S(\OO_E)(y_2)\not\simeq S(\PP^d_x)(y_1)$$
and
$S(\OO_E)(y_2)\not\simeq S(\OO_E)(y_1)$.
Since a simple vector bundle is determined up to an isomorphism
by its determinant, it suffices to check that
$$\det(S(\OO_E))(d(y_2-y_1))\not\simeq\det(S(\PP^d_x)).$$
$$\det(S(\OO_E))(d(y_2-y_1))\not\simeq\det(S(\OO_E)).$$
It is easy to see that we have an equality $y_2-y_1=\pm dy$ in the group
$E(k)$. Changing $S$ by $[-id_E]^*S$ if necessary we can assume that
$y_2-y_1=dy$. Then considering the action of $S$ on $K_0(E)$ we derive
the isomorphism 
$$\det(S(\PP^d_x))\simeq \det(S(\OO_E))(x'-e)$$
where $x'=dx$ in $E(k)$. Our assertion follows.
\ed

Thus, in the case $k=\C$ using some uniformization $\pi:\C\ra E$ we
can consider the functions
$$r_V(u,v):=r_V(\pi(u),\pi(v))$$
and
$$r_V(u):=r_V(\pi(u))$$
as meromorphic solutions of the AYBE and CYBE
respectively satisfying some additional conditions
(namely, the unitarity and the non-degeneracy conditions).

In particular, $r_V(u)$ is a solution of CYBE
satisfying all the additional conditions imposed by Belavin and
Drinfeld in \cite{BD}.
The explicit formulas of section \ref{KWsec}
imply that $r_V(v)$ has poles at the points of the
lattice $\pi^{-1}(E_d)$ (and is periodic with respect to the lattice
$\pi^{-1}(0)$). In order to find the place of $r_V(u)$ in
Belavin-Drinfeld classification we have to determine the automorphisms
$$A_{\ga}:\pgl(H^0(E,V))\ra\pgl(H^0(E,V))$$
for all $\ga\in \pi^{-1}(E_d)$ such that
$$r_V(u+\ga)=(A_{\ga}\otimes 1)r_V(u)$$
(see Prop.4.3 of \cite{BD}). Note that by periodicity of $r_V(u)$
with respect to $\pi^{-1}(0)$ the automorphism $A_{\ga}$ depends
only on $\pi(\ga)\in E_d$.

Let $H$ be the Heisenberg group associated with $V$. Recall that
$H$ is the central extension of $E_d$ (the subgroup of points
of order $d$ in $E$) by $\G_m$. Points of $H$ are pairs
$(x,\a)$ where $x\in E_d$, $\a:V\ra t_x^*V$ is an isomorphism.
The space $H^0(E,V)$ is an irreducible representation of $H$ in a
natural way. This induces a natural action of $E_d=H/\G_m$ on 
$\pgl(H^0(E,V))$. It is easy to see that the automorphism
$A_{\ga}$ above is given by the action of $\pi(\ga)\in E_d$.

The solution $r_V(u)$ gets replaced by an equivalent one
if we replace $V$ by $T(V)$ where $V$ is any autoequivalence of
$D^b(E)$ preserving $\OO_E$. Thus, the only data on which $r_V(u)$
depends are $(d=\deg(V),r=\rk(V)\mod d)$. Note that the rank $r$ is
relatively prime to $d$ since $V$ is simple. It follows
that the solutions for $\pgl_d$ are numbered by $(\Z/d\Z)^*$.
The choice of $r\in (\Z/d\Z)^*$ precisely corresponds to a choice
of a primitive $d$-th root of unity in Belavin-Drinfeld's picture.

\subsection{Explicit formulas}
\label{KWsec}

Now we assume that $k=\C$ and write explicit formulas
for the above solutions.
The elliptic solutions of the AYBE can be expressed
in terms of the Kronecker function
\begin{equation}\label{Kr-id}
F(u,v)=\frac{\th'_{11}(0)}{2\pi i}\cdot\frac{\th_{11}(u+v)}
{\th_{11}(u)\th_{11}(v)}
\end{equation}
where
$$\th_{11}(u,\tau)=\sum_{n\in\Z}(-1)^n\exp(\pi i (n+\frac{1}{2})^2\tau+
2\pi i (n+\frac{1}{2})u),$$
$\th'_{11}$ is the derivative of $\th_{11}(u,\tau)$ with respect to $u$.
When we want to stress the dependance of $F$ on $\tau$ we will write
$F(u,v,\tau)$. Kronecker discovered the following series expansion:
$$F(u,v)=-\sum_{(m+\frac{1}{2})(n+\frac{1}{2})>0}
\sign(m+\frac{1}{2})\exp(2\pi i (mn\tau+mv+nu))$$
where $m,n$ are integers, $0<\Im(u),\Im(v)<\Im(\tau)$.
Let us introduce a little bit more notation. For a pair of rational numbers
$(p,q)$ we set
\begin{equation}\label{Kr-char}
F_{p,q}(u,v)=\exp(2\pi i (pq\tau+pv+qu))F(u+p\tau,v+q\tau).
\end{equation}
For $0<\Im(u),\Im(v)<\eps$ where $\eps$ is sufficiently small, one has
$$F_{p,q}(u,v)=-\sum_{(m,n)\in\Z^2+(p,q),(m+\eps)(n+\eps)>0}\sign(m+\eps)
\exp(2\pi i (mn\tau+mv+nu)).$$
Note that we have the symmetry relation
$$F_{p,q}(u,v)=F_{q,p}(v,u).$$
This kind of series appear in the computation of triple Fukaya compositions
corresponding to the Massey products defining $r_V(u,v)$.

Let us consider first the case $r=1$, so $V=L$ is a line bundle of degree $d$.
We denote by $(e_i,i\in\Z/d\Z)$ the natural basis in $H^0(E,L)$
consisting of theta-functions with characteristics.
Let $e^*_i$ be the dual basis in $H^0(E,L)^*$.
Then using the correspondence
between our Massey products and triple Fukaya compositions (see \cite{P-Mas})
one can derive the following formula: 
$$m_3(e_i,e^*_j,e_k)=F_{\frac{i-j}{d},\frac{j-k}{d}}(du,-dv,d\tau)e_{i-j+k}.$$
Hence,
\begin{equation}\label{AYBEfor1}
r_L(u,v)=\sum_{j-i=i'-j'}F_{\frac{j-i}{d},\frac{i-j'}{d}}(du,-dv,d\tau)
e_{ij}\otimes e_{i'j'}
\end{equation}
where $e_{ij}$ is the standard basis in the matrix algebra $\Mat(d,\C)$.
In the simplest case when $d=1$ we obtain just the function $F(u,-v)$, so
the AYBE in this case specializes to the following identity:
\begin{equation}
F(-u',v)F(u+u',v+v')-F(u+u',v')F(u,v)+
F(u,v+v')F(u',v')=0.
\end{equation}

To find formulas for the corresponding solutions of the CYBE
we project the tensor $r_L(u,v)\in\Mat(d,\C)\otimes\Mat(d,\C)$ to
$\splin_d\otimes\splin_d$ and then set $u=0$.
Using the above formula for $r_L(u,v)$ we obtain
\begin{equation}\label{prfor}
\ov{r}(v):=(\pr\otimes\pr)(r_L(u,v))=\sum_{j-i=i'-j'\neq 0}
F_{\frac{j-i}{d},\frac{i-j'}{d}}(du,-dv,d\tau)
e_{ij}\otimes e_{i'j'}+\sum_{i,i'} G_{i-i'}(du,-dv,d\tau)
e_{ii}\otimes e_{i'i'}
\end{equation}
where 
$$G_{j}(x,y,\tau)=F_{0,\frac{j}{d}}(x,y,\tau)-\frac{1}{d}\cdot
\sum_{k\in\Z/d\Z}F_{0,\frac{k}{d}}(x,y,\tau)$$

When passing to the limit $u\ra 0$ in the formula (\ref{prfor})
we are going to use the following relation between the Kronecker
function $F(u,v)$ and the Weierstrass zeta-function observed in
\cite{P-KW}. Let $\zeta(x)=\zeta(x,\tau)$
denotes the Weierstrass zeta-function
associated with the lattice $\Z+\Z\tau$.
Then according to \cite{P-KW}, Cor.1.2, we have
$$\left(2\pi i F(x,y)-\frac{1}{x}\right)|_{x=0}=\zeta(y)-y\eta_1$$
where $\eta_1=2\zeta(\frac{1}{2})$.
It follows that for any function $g:\Z/d\Z\ra\C$ with
$\sum_{j\in\Z/d\Z}g(j)=0$ we have
$$\left(2\pi i \sum_{j\in\Z/d\Z}g(j)F_{0,\frac{j}{d}}(x,y)\right)|_{x=0}=
\sum_{j\in\Z/d\Z}g(j)\zeta(y+\frac{j}{d}\tau)+
(2\pi i-\eta_1\tau)\sum_{j\in\Z/d\Z}\frac{g(j)j}{d}.$$
Using the Legendre relation $\eta_1\tau-\eta_2=2\pi i$,
where $\eta_2=2\zeta(\frac{\tau}{2})$, we can rewrite this formula as
follows:
\begin{equation}
\left(2\pi i \sum_{j\in\Z/d\Z}g(j)F_{0,\frac{j}{d}}(x,y)\right)|_{x=0}=
\sum_{j\in\Z/d\Z}g(j)\zeta_{0,\frac{j}{d}}(y,\tau).
\end{equation}
where we use the notation (\ref{zeta-char}). In particular, we obtain
$$2\pi i G_{j}(0,y,\tau)=\zeta_{0,\frac{j}{d}}(y,\tau)-\frac{1}{d}\cdot
\sum_{k\in\Z/d\Z}\zeta_{0,\frac{k}{d}}(y,\tau).$$
Now the expression for the solutions of the CYBE
takes form
\begin{equation}\label{CYBEfor}
\begin{array}{l}
2\pi i \ov{r}(v)=\sum_{j-i=i'-j'\neq 0}
2\pi i F_{\frac{j-i}{d},\frac{i-j'}{d}}(0,-dv,d\tau)
e_{ij}\otimes e_{i'j'}+\\
\sum_{i,i'} (\zeta_{0,\frac{i-i'}{d}}(-dv,d\tau)-\frac{1}{d}
\sum_{k\in\Z/d\Z}\zeta_{0,\frac{k}{d}}(-dv,d\tau))
e_{ii}\otimes e_{i'i'}.
\end{array}
\end{equation}
Using formulas (\ref{formula1}) and (\ref{formula2})
we can rewrite this as follows:
\begin{equation}\label{CYBEforbis}
\begin{array}{l}
2\pi i \ov{r}(v)=\sum_{j-i=i'-j'\neq 0}
\sum_{a\in\Z/d\Z}\exp(-2\pi i\frac{a(j-i)}{d})
[\zeta_{\frac{a}{d},\frac{i-j'}{d}}(-v,\tau)-
\zeta_{\frac{a}{d},0}(\frac{i-j}{d}\tau,\tau)]
e_{ij}\otimes e_{i'j'}+\\
\sum_{i,i'}[\frac{1}{d}
 \sum_{a\in\Z/d\Z}\zeta_{\frac{a}{d},\frac{i-i'}{d}}(-v,\tau)-
\frac{1}{d^2}\sum_{a,b\in\Z/d\Z}\zeta_{\frac{a}{d},\frac{b}{d}}(-v,\tau)]
e_{ii}\otimes e_{i'i'}.
\end{array}
\end{equation}

The case $r>1$ can be easily reduced to the case $r=1$ using a representation
of the bundle $V$ as the direct image of a line bundle $L$ under the isogeny
$\C/\Z+r\tau\Z\ra\C/\Z+\tau\Z$. It is easy to see that in this situation
one has
$$r_V(u,v;\tau)=r_L(ru,v,r\tau).$$

%

\section{Trigonometric solutions for $\splin_2$}
\label{trig-sec}

It turns out that computations of Massey products are easier in the
case of a reducible curve. Also in order to obtain all non-degenerate
solutions of the CYBE for $\splin_2$ it is necessary to consider a curve with
$2$ components. Because of this we chose to study the solutions of
the AYBE and the CYBE arising from simple bundles of rank $2$ on such a
curve.

\subsection{Construction of simple bundles of rank $2$ on a reducible curve}

Let $C=C_1\cup C_2$ be the union of two $\P^1$'s glued (transversally) by
two points. In other words, $C_1=C_2=\P^1$ and the point $0$ (resp.
$\infty$) on $C_1$ is identified with the point $0$ (resp. $\infty$)
on $C_2$. A vector bundle $V$ on $C$ is given by the following data:
$$(V_1,V_2,\a_0:V_{1,0}\wt{\ra}V_{2,0},
\a_{\infty}:V_{1,\infty}\wt{\ra} V_{2,\infty})$$
where $V_i$ is a bundle on $C_i$, $i=1,2$, $V_{i,x}$ denotes the fiber
of $V_i$ at the point $x$. For each $\la\in k^*$ let us define the rank-$2$
bundle $V^{\la}$ on $C$ as follows:
$$V^{\la}_1=\OO_{\P^1}\oplus\OO_{\P^1},$$
$$V^{\la}_2=\OO_{\P^1}\oplus\OO_{\P^1}(1),$$
$$\a_0=\id,\ \a_{\infty}=
S_{\la}:=\left(\matrix 0 & \la\\ 1 & 0\endmatrix\right).$$
Here we use the trivialization of $\OO_{\P^1}(1)$ at $0$
(resp. $\infty$) induced by the standard trivialization of $\OO_{\P^1}(1)$
on the complement to $\infty$ (resp. $0$). 

\begin{lem} The bundle $V^{\la}$ is simple.
\end{lem}

\Pf . An endomorphism of $V^{\la}$ is given by a pair of endomorphisms
$f_1:V^{\la}_1\ra V^{\la}_1$ and $f_2:V^{\la}_2\ra V^{\la}_2$,
such that $f_1(0)=f_2(0)$ (this follows from $\a_0=\id$) and
$$f_2(\infty)S_{\la}=S_{\la}f_1(\infty).$$
Note that $f_1$ has constant coefficients so $f_1=f_1(0)=f_1(\infty)$.
The endomorphism $f_2$ is lower-triangular (since 
$\Hom(\OO_{\P^1}(1),\OO_{\P^1})=0$), 
hence, $f_1$ and $S_{\la}f_1 S_{\la}^{-1}$ are both lower-triangular
which implies that $f_1$ is diagonal. Notice that the
diagonal part of $f_2$ is constant, so we deduce that
$$f_1=S_{\la}f_1 S_{\la}^{-1}$$
which is possible only if $f_1$ is proportional to the identity.
Finally, it is easy to see that $f_2$ is completely determined by
$f_2(0)$ and $f_2(\infty)$, so the only endomorphisms of $V$ are scalar
multiples of the identity.
\ed

\subsection{Computation}

Now we are going to apply theorem \ref{bundlethm} to compute the solutions
of the AYBE and the CYBE associated with bundles $V^{\la}$.
For this we have to describe the space of morphisms
$\Hom(V^{\la_1},V^{\la_2}(y))$, where $\la_i\in k^*$, $y$ is a smooth
point of $C$. There are two different cases to consider depending
on whether $y\in C_1$ or $y\in C_2$.

\noindent {\bf Case 1}. $y\in C_1$.
Then a morphism $V^{\la_1}\ra V^{\la_2}(y)$ is given by a pair of
morphisms on $\P^1$:
$$A:\OO_{\P^1}\oplus\OO_{\P^1}\ra \OO_{\P^1}(y)\oplus\OO_{\P^1}(y),$$
$$B:\OO_{\P^1}\oplus\OO_{\P^1}(1)\ra \OO_{\P^1}\oplus\OO_{\P^1}(1)$$
satisfying the conditions $A_0=B_0$ and
$$S_{\la_2}A_{\infty}=B_{\infty}S_{\la_1}.$$
We claim that such a morphism is
completely determined by $B$ which can be arbitrary.
Indeed, considering $A$ as an endomorphism of $\OO_{\P^1}^2$
with a pole of the first order at $y$ we can write it uniquely in the form
$$A=\frac{1}{z-y}\cdot A'+\frac{z}{z-y}\cdot A''$$
where $A',A''$ are some regular endomorphisms of $\OO_{\P^1}^2$,
$z=\frac{z_1}{z_0}$. Now we have
$$A_0=-\frac{A'}{y}, \ A_{\infty}=A'',$$
hence $A$ is uniquely recovered from $A_0$ and $A_{\infty}$.
Thus, to every $B\in\End(\OO_{\PP^1}\oplus\OO_{\P^1}(1))$
we can associate the morphism $(A,B):V^{\la^1}\ra V^{\la_2}(y)$
with
$$A=\frac{y}{y-z}\cdot B_0+\frac{z}{z-y}\cdot S_{\la_2}^{-1}B_{\infty}
S_{\la_1}.$$
In this description the residue morphism
$$\Res_y:\Hom(V^{\la_1},V^{\la_2}(y))\ra\Mat(2,k)$$
is given by the formula
$$B\mapsto S_{\la_2}^{-1}B_{\infty}S_{\la_1}-B_0$$
(here we use a local trivialization of $\om_C$ given by the form
$\frac{dz}{z}$). Let us write
$$B=\left(\matrix a & 0 \\ bz_0+cz_1 & d\endmatrix\right).$$
Then we have
$$\Res_y:B\mapsto 
\left(\matrix d-a & \la_1 c\\ -b & \la_1\la_2^{-1}a-d\endmatrix\right).$$
On the other hand, if $y_1,y_2\in C_1$ are distinct points then
after applying the above computation to $y=y_1$
we can consider the evaluation map 
$$\ev_{y_2}:\Hom(V^{\la_1},V^{\la_2}(y_1))\ra\Mat(2,k):
B\mapsto \frac{y_1}{y_1-y_2}\cdot B_0+\frac{y_2}{y_2-y_1}\cdot
S_{\la_2}^{-1}B_{\infty}S_{\la_1}.$$
Thus, we can compute the map
\begin{align*}
&\ev_{y_2}\circ\Res_{y_1}^{-1}:\Mat(2,k)\ra\Mat(2,k):
\left(\matrix a & b\\ c & d\endmatrix\right)\mapsto\\
&\frac{y_1}{y_1-y_2}\cdot
\left(\matrix \frac{a+d}{\la_1\la_2^{-1}-1} & 0 \\ -c &
\frac{\la_1\la_2^{-1}a+d}{\la_1\la_2^{-1}-1}\endmatrix\right)+
\frac{y_2}{y_2-y_1}\cdot
\left(\matrix \frac{\la_1\la_2^{-1}a+d}{\la_1\la_2^{-1}-1} & b\\ 0 &
\la_1\la_2^{-1}\cdot\frac{a+d}{\la_1\la_2^{-1}-1}\endmatrix\right).
\end{align*}
Note that this map depends only on $\la=\la_1\la_2^{-1}$ and
$\mu=y_1y_2^{-1}$. Thus, from theorem \ref{bundlethm} we obtain
the following solution of the AYBE (where $\la$ and $\mu$ should
be considered as multiplicative variables which are exponents
of the additive variables appearing in (\ref{AYBE})):
\begin{equation}\label{AYBEsl2-1}
\begin{array}{l}
r(\la,\mu)=\frac{1}{(1-\la)(1-\mu)}
\left((\mu e_{11}-e_{22})\otimes(e_{11}+\la e_{22})+
(-\la e_{11}+\mu e_{22})\otimes(e_{11}+e_{22})\right)+\\
\frac{1}{1-\mu}e_{21}\otimes e_{12}+\frac{\mu}{1-\mu}e_{12}\otimes e_{21}.
\end{array}
\end{equation}

Projecting this tensor to $\splin_2$ we obtain the corresponding solution
of the CYBE:
\begin{equation}\label{CYBEsl2-1}
r(\mu)=\frac{1+\mu}{4(1-\mu)}h\otimes h+\frac{e_{21}\otimes e_{12}+\mu
e_{12}\otimes e_{21}}{1-\mu}.
\end{equation}
where $h=e_{11}-e_{22}$.

\noindent {\bf Case 2}. $y\in C_2$.
Then a morphism $V^{\la_1}\ra V^{\la_2}(y)$ is given by a pair of
morphisms on $\P^1$:
$$A:\OO_{\P^1}\oplus\OO_{\P^1}\ra \OO_{\P^1}\oplus\OO_{\P^1},$$
$$B:\OO_{\P^1}\oplus\OO_{\P^1}(1)\ra (\OO_{\P^1}\oplus\OO_{\P^1}(1))(y)$$
satisfying the conditions $A=B_0$ and $S_{\la_2}A=B_{\infty}S_{\la_1}.$
Such a morphism is completely determined by $B$ which should satisfy
the condition 
$$B_0=S_{\la_2}^{-1}B_{\infty}S_{\la_1}.$$
Considering $B$ as an endomorphism of $\OO_{\P^1}\oplus\OO_{\P^1}(1)$
with a pole of the first order at $y$ we can write it in the form
$$B=\frac{1}{z-y}\cdot B'+\frac{z}{z-y}\cdot B''+
\left(\matrix 0 & \frac{t}{z_1-z_0y} \\ 0 & 0\endmatrix\right)$$
where $B'$ and $B''$ are regular endomorphisms of $\OO_{\P^1}\oplus
\OO_{\P^1}(1)$, $t\in k$, $\frac{1}{z_1-z_0y}$ is a section of 
$\OO_{\P^1}(-1)$ with the pole at $y$. 
However, this presentation is non-unique:
we can add to $B'$ a lower-triangular endomorphism vanishing at $0$
and change $B''$ appropriately. To get rid of this ambiguity we
impose the condition that $B'_{\infty}$ is diagonal. Then $B'$ and
$B''$ are unique. Furthermore, in this case $B'$ is uniquely
determined by $B'_0$. On the other hand, we have
$$B_0=-\frac{B'_0+te_{12}}{y},$$ 
$$B_{\infty}=B''_{\infty}+te_{12},$$ 
hence we get the equation
$$B'_0+te_{12}=-yS_{\la_2}^{-1}(B''_{\infty}+te_{12})S_{\la_1}.$$
Solving this equation for $B'_0$ and $t$ we obtain that 
for
$$B''=\left(\matrix a'' & 0 \\ b''z_0+c''z_1 & d''\endmatrix\right)$$
one has $t=-y\la_1c''$ and 
$$B'_0=-y\cdot\left(\matrix d'' & 0 \\ -y\la_1\la_2^{-1}c'' & 
\la_1\la_2^{-1}a'' \endmatrix\right).$$ 
Thus, all the data can be recovered from $B''$ which can be arbitrary.
Now we can compute the map $\Res_y$. Notice that the difference
from the previous case is that we have to choose a trivialization
of $V^{\la_1}$ and $V^{\la_2}$ at $y$ (since now $y$ belongs to the
component $C_2$ on which these bundles are non-trivial).
Our choice for $V^{\la}$ will correspond to the trivialization of 
$\OO_{\P^1}(1)$ at $y$ given by the non-vanishing section
$f_{\la}^{-1}z_0$, where $f_{\la}$ is some invertible function on 
$\P^{1}-\{0,\infty\}$. 
Then using $B''$ as a coordinate on $\Hom(V^{\la_1},V^{\la_2})$
we obtain
$$\Res_y(B'')=\frac{B'_y}{y}+B''_y+\frac{t}{y}e_{12}=
\left(\matrix a''-d'' & -f_{\la_1}^{-1}(y)\la_1c'' \\
f_{\la_2}(y)(b''+y(1+\la_1\la_2^{-1})c'') & 
d''-\la_1\la_2^{-1}a''\endmatrix\right)$$
On the other hand, using the above construction for $y=y_1$ and taking
a point $y_2\neq y_1$ in $C_2$ we can compute the evaluation map 
\begin{align*}
&\ev_{y_2}:\Hom(V^{\la_1},V^{\la_2}(y_1))\ra\Mat(2,k):\\
&B''\mapsto \frac{1}{y_2-y_1}B'_{y_2}+\frac{y_2}{y_2-y_1}B''_{y_2}+
\frac{t}{y_2-y_1}e_{12}=\\
&\frac{y_2}{y_2-y_1}\cdot
\left(\matrix a''-y_1y_2^{-1}d'' & -f_{\la_1}^{-1}(y_2)\la_1c'' \\
f_{\la_2}(y_2)(b''+(y_2+y_1^2y_2^{-1}\la_1\la_2^{-1})c'') &
d''-y_1y_2^{-1}\la_1\la_2^{-1}a''\endmatrix\right).
\end{align*}
Finally, we compute the map
\begin{align*}
&\ev_{y_2}\circ\Res_{y_1}^{-1}:\Mat(2,k)\ra\Mat(2,k):
\left(\matrix a & b\\ c & d\endmatrix\right)\mapsto\\
&\frac{1}{1-\mu}\cdot
\left(\matrix \frac{(1-\mu\la)a+(1-\mu) d}
{1-\la} & f_{\la_1}(y_2)^{-1}f_{\la_1}(y_1)b\\ 
f_{\la_2}(y_2)[f_{\la_1}(y_1)(y_1-y_2)
(1-\mu\la)\la_1^{-1}b+f_{\la_2}(y_1)^{-1}c] &
\frac{(1-\mu)\la a+(1-\mu\la)d}
{1-\la}\endmatrix\right)
\end{align*}
where we denoted $\la=\la_1\la_2^{-1}$, $\mu=y_1y_2^{-1}$.
Now we observe that if we set
$$f_{\la}(y)=\la^{\frac{1}{2}}y^{-\frac{1}{2}}$$
then the above matrix will depend only on $\la$ and $\mu$.
Thus, we obtain the following solution of the AYBE (in the
multiplicative notation):
\begin{equation}\label{AYBEsl2-2}
\begin{array}{l}
r(\la,\mu)=\frac{1-\la\mu}{(1-\la)(1-\mu)}
(e_{11}\otimes e_{11}+e_{22}\otimes e_{22})
+\frac{\la e_{11}\otimes e_{22}+e_{22}\otimes e_{11}}{1-\la}+\\
\frac{\mu^{-\frac{1}{2}}}{1-\mu} e_{21}\otimes e_{12}+
\frac{\mu^{\frac{1}{2}}}{1-\mu} e_{12}\otimes e_{21}+
((\la\mu)^{\frac{1}{2}}-(\la\mu)^{-\frac{1}{2}}) e_{21}\otimes e_{21}.
\end{array}
\end{equation}

Applying the projection to $\splin_2$ and setting $\la=1$ we get
the following solution of the CYBE:
\begin{equation}\label{CYBEsl2-2}
r(\mu)=\frac{1+\mu}{4(1-\mu)}h\otimes h+
\frac{\mu^{-\frac{1}{2}}e_{21}\otimes e_{12}+\mu^{\frac{1}{2}}
e_{12}\otimes e_{21}}{1-\mu}+ (\mu^{\frac{1}{2}}-\mu^{-\frac{1}{2}})
e_{21}\otimes e_{21}.
\end{equation}
where $h=e_{11}-e_{22}$.
It is easy to see that our solutions (\ref{CYBEsl2-1}) and
(\ref{CYBEsl2-2}) are equivalent to the solutions (6.9) and (6.10) in 
\cite{BD} which represent two distinct equivalence classes
of non-degenerate trigonometric solutions
of the CYBE for $\splin_2$. Note that we actually constructed a
solution $r_{y_1,y_2}$ of the equation (\ref{CYBE}) depending on parameters 
$y_1,y_2\in C^{reg}$ (with a pole at $y_1=y_2$) which specializes to the above
two solutions when $y_i$ vary in one of the two components of $C$.

\section{Scalar solutions of AYBE}
\label{scalar}

In this section we are going to study the 
equation (\ref{AYBE}) in the case when $n=1$,
i.e. when $r(u,v)$ is $\C$-valued. 

\begin{thm}\label{scthm} Let $r(u,v)$ be a non-zero 
meromorphic function in the neighborhood of
$(0,0)$ satisfying the equations
$$r(-u',v)r(u+u',v+v')-r(u+u',v')r(u,v)+
r(u,v+v')r(u',v')=0,$$
$$r(-u,-v)=-r(u,v).$$
Then there exist constants $c_1,c_3,c_4\in\C^*$ and $c_2\in\C$ such
that $c_1\exp(c_2uv)r(c_3u,c_4v)$ is one of the following functions:

\noindent 1) $F_{\tau}(u,v)=F(u,v,\tau)$ (Kronecker's function),\newline

\noindent 2) $F_{\infty}(u,v):=\frac{\exp(v)-\exp(u)}{(\exp(u)-1)(\exp(v)-1)}$,
\newline

\noindent 3) $\frac{a}{u}+\frac{b}{v}$,\ $a,b\in\C$. 
\end{thm}

\Pf . Assume first the divisor of poles of $r$ doesn't contain $u=0$.
Substituting $u=0$ in the equation we obtain
\begin{equation}\label{aux1}
r(-u',v)r(u',v+v')-r(u',v')r(0,v)+r(0,v+v')r(u',v')=0.
\end{equation}
Substiting $u'=0$ we get
$$(r(0,v)+r(0,v'))r(0,v+v')=r(0,v)r(0,v').$$
Note that $r(0,v)$ is not identically zero: otherwise
(\ref{aux1}) would imply that $r(u,v)$ is identically zero.
Hence, we can write the last equation as
$$\frac{1}{r(0,v+v')}=\frac{1}{r(0,v)}+\frac{1}{r(0,v')}.$$
Thus, multiplying $r$ by a constant we can assume that
$$r(0,v)=\frac{1}{v}.$$
Substituting this in (\ref{aux1}) we obtain
$$r(-u',v)r(u',v+v')=r(u',v')(\frac{1}{v}-\frac{1}{v+v'}).$$
Using the equality $r(-u',v)=-r(u',-v)$ we can rewrite this
as follows
$$-r(u',-v)r(u',v+v')v(v+v')=r(u',v')v'.$$
This implies that 
\begin{equation}\label{aux2}
r(u,v)v=\exp(c(u)v)
\end{equation}
for some meromorphic function $c(u)$. 
Substituting this in the original equation we get
$$\frac{\exp(c(-u')v+c(u+u')(v+v'))}{v(v+v')}-
\frac{\exp(c(u+u')v'+c(u)v)}{vv'}+
\frac{\exp(c(u)(v+v')+c(u')v')}{v'(v+v')}=0
$$
Multiplying by $v+v'$ and collecting terms with $1/v$ and $1/v'$ we
get
$$\frac{\exp((c(u+u')-c(u)-c(u'))v)-1}{v}=
\frac{1-\exp((c(u)+c(u')-c(u+u'))v')}{v'}.$$
This immediately implies that
$$c(u+u')=c(u)+c(u'),$$
hence 
$$r(u,v)=\frac{\exp(cuv)}{v}$$
for some constant $c$ which leads to case 3).

Now let us assume that $r$ has pole along $u=0$ of order $k>0$.
Writing $r$ in the form 
$r(u,v)=\sum_{i\le -k}r_i(v)u^i$ and substituting in the equation
we obtain that
$$r_{-k}(v)r_{-k}(v+v')(-u')^{-k}(u+u')^{-k}-
r_{-k}(v')r_{-k}(v)u^{-k}(u+u')^{-k}+
r_{-k}(v+v')r_{-k}(v')u^{-k}(u')^{-k}=0.$$
It is easy to see that this is possible only if $k=1$ and $r_{-1}(v)$
is constant. 
Multiplying $r$ be a constant we can assume that
$$r(u,v)=\frac{1}{u}+r_0(v)+r_1(v)u+r_2(v)u^2+\ldots$$
Note that similar arguments work for $v$ instead of $u$,
so we can assume that $r_0(v)$ has pole of order $1$ are zero.

Now we claim that the terms $r_i$ with $i\ge 2$ are uniquely
determined by $r_0$ and $r_1$. Indeed, let us check that
the term $r_n$ for $n\ge 2$ can be recovered from the previous term.
Collecting terms of the main equation which have total degree $n-1$ in
$u$ and $u'$ we get 
$$r_n(v)[\frac{(-u')^n}{u+u'}-\frac{u^n}{u+u'}]+
r_n(v')[-\frac{(u+u')^n}{u}+\frac{(u')^n}{u}]+
r_n(v+v')[-\frac{(u+u')^n}{u'}+\frac{u^n}{u'}]=\ldots
$$
where the RHS contains only $r_i$ with $i\le n-1$.
It is easy to check that if $n\ge 3$ then the polynomials in $u,u'$
$$\frac{(-u')^n-u^n}{u+u'};\ 
\frac{-(u+u')^n+(u')^n}{u};\
\frac{-(u+u')^n+u^n}{u'}$$
are linearly independent (e.g. one can check this by looking at coefficients
with $u^{n-1}$, $u^{n-2}u'$ and $(u')^{n-1}$). Therefore, for $n\ge 3$
the term $r_n$ is recovered from the previous terms. For $n=2$ the above
equation takes form
$$-u(r_2(v)+r_2(v'))+u'(r_2(v)-r_2(v+v'))=\ldots,$$
hence, $r_2$ is uniquely recovered from $r_0$ and $r_1$.
For $n=1$ we get the following relation
$$r(-u',v)r(u+u',v+v')-r(u+u',v')r(u,v)+
r(u,v+v')r(u',v')=0,$$
\begin{equation}\label{aux3}
r_0(v)r_0(v+v')-r_0(v')r_0(v)+r_0(v+v')r_0(v')=
r_1(v)+r_1(v')+r_1(v+v')
\end{equation}
Using the rescaling of the form 
$$r(u,v)\mapsto c\cdot\exp(c'uv)r(cu,c''v)$$
we can achieve rescaling of $r_0$ of the form
$$r_0(v)\mapsto cr_0(c''v)+c'v.$$
Thus, we can assume that the Laurent expansion of $r_0$ at $0$ has form
\begin{equation}\label{form}
r_0(v)=\frac{1}{v}+c_3v^3+c_5v^5+\ldots
\end{equation}
where $c_3$ is equal to $1$ or $0$ (recall that $r_0$ is odd).
Note that the LHS in (\ref{aux3}) doesn't have pole at $v=0$.
Hence, $r_1$ is regular at $0$ and taking the limit of (\ref{aux3})
as $v\ra 0$ we get
$$r'_0(v')+r_0(v')^2=r_1(0)+2r_1(v').$$
Using the Laurent expansion of $r_0$ at $0$ we see that
the LHS of this equality tends to zero as $v'\ra 0$. Hence,
$r_1(0)=0$ and we get
$$r_1(v')=\frac{1}{2}(r'_0(v')+r_0(v')^2).$$
In particular, $r_1$ is determined by $r_0$. Substituting this
expression for $r_1$ into (\ref{aux3}) we obtain the following
functional equation on $r_0$:
$$2r_0(v)r_0(v+v')-2r_0(v')r_0(v)+2r_0(v+v')r_0(v')=
r'_0(v)+r'_0(v')+r'_0(v+v')+r_0(v)^2+r_0(v')^2+r_0(v+v')^2,$$
which can be rewritten as
\begin{equation}\label{aux4}
(r_0(v)+r_0(v')-r_0(v+v'))^2+r'_0(v)+r'_0(v')+r'_0(v+v')=0.
\end{equation}
We are looking for solutions of this equation which are meromorphic in
the neighborhood of zero and have form (\ref{form}).\footnote{
Meromorphic solutions of (\ref{aux4}) were described by L.~Carlitz in 
\cite{C}. For completeness we give an independent argument.} Substituting
the expansion (\ref{form}) in the equation one can easily see that
any solution is uniquely determined by the coefficients $(c_3,c_5)$.
The rescaling $r_0(v)\mapsto cr_0(cv)$ for $c\in\C^*$ leads to the rescaling 
$(c_3,c_5)\mapsto (c^4c_3, c^6c_5)$. Note that there is a (unique) solution
with $c_3=c_5=0$, namely, $r_0(v)=\frac{1}{v}$ (this corresponds to
$r(u,v)=\frac{1}{u}+\frac{1}{v}$), so from now on we will assume
that $(c_3,c_5)\neq (0,0)$. Then up to rescaling a solution $r_0$ is 
characterized by the parameter 
$$C(r_0)=\frac{c_5^2}{c_3^3}$$ 
which takes values in $\C\cup\infty$. Now we claim that from solutions
$r(u,v)=2\pi i F_{\tau}(u,v)$ one gets all values of $C(r_0)$ except for
$-\frac{20}{49}$ while from the trigonometric solution 
$r(u,v)=F_{\infty}(u,v)$ one gets the exceptional value $-\frac{20}{49}$.
Indeed, the Laurent expansion $2\pi i F_{\tau}(u,v)$ has form
$$2\pi i F_{\tau}(u,v)=\frac{1}{u}+[\frac{1}{v}-2G_2(\tau)v
-G_4(\tau)\frac{v^3}{3}-G_6(\tau)\frac{v^5}{60}+\ldots]+\ldots$$
where 
$$G_{k}=-\frac{B_k}{2k}+\sum_{m,n\ge 1}m^{k-1}q^{mn}$$
are the Eisenstein series (here $q=\exp(2\pi i\tau)$).
Thus, for this solution we have
$$C(r_0)=-\frac{27G_6(\tau)^2}{60^2G_4(\tau)^3}.$$
Recall that the $j$-invariant is defined by the formula
$$j(\tau)=\frac{g_2(\tau)^3}{g^2(\tau)^3-27g_3(\tau)^2}$$
where 
$$g_2(\tau)=60\sum_{(m,n)\neq(0,0)}\frac{1}{(m\tau+n)^4},$$
$$g_3(\tau)=140\sum_{(m,n)\neq(0,0)}\frac{1}{(m\tau+n)^6}.$$
We have the following relations:
$$G_4(\tau)=\frac{g_2(\tau)}{20(2\pi)^4},$$
$$G_6(\tau)=-\frac{3g_3(\tau)}{7(2\pi )^6}.$$
It follows that
$$C(r_0)=-\frac{20}{49}(1-j(\tau)^{-1}).$$
Since $j(\tau)$ takes all complex values (including $0$), we obtain
all values of the parameter $C(r_0)$ (including $\infty$) except for 
$-\frac{20}{49}$. Finally, for the solution $r(u,v)=F_{\infty}(u,v)$
we obtain
$$r_0(v)=\frac{\coth(v/2)}{2}=\sum_{n\ge 0}\frac{B_n}{n!}v^{n-1},$$
hence
$$C(r_0)=(\frac{B_6}{6!})^2\cdot(\frac{4!}{B_4})^3=-\frac{20}{49}.$$
\ed

\section{Reconstructing solutions of AYBE from solutions of CYBE}
\label{reconstr-sec}

Recall that according to Lemma \ref{limit} if $r(u,v)$ is a unitary
solution of the AYBE 
then the limit $\ov{r}(v):=(\pr\otimes\pr)(r(u,v))|_{u=0}$ (if exists) is
a solution of the CYBE with values in $\splin_n$. 
In this section we study the question to
which extent $r(u,v)$ is determined by $\ov{r}(v)$.

\begin{thm}\label{reconstr} 
Consider unitary solutions of the AYBE with values in $\Mat_n(\C)$
which have Laurent expansion near $u=0$ of the form
$$r(u,v)=\frac{1\otimes 1}{u}+r_0(v)+r_1(v)u+\ldots.$$
Assume that the corresponding solution 
$\ov{r}(v):=(\pr\otimes\pr)(r_0(v))$ of the CYBE 
has no infinitesimal symmetries and
that the tensor $\ov{r}(v)$ has rank $>2$ for generic $v$.
Then $r(u,v)$ can be uniquely recovered from $\ov{r}(v)$ up to a rescaling
$r(u,v)\mapsto \exp(cuv)r(u,v)$, where $c\in\C$.
In other words, two unitary solutions of the AYBE in the above form
differ by a factor of the form $\exp(cuv)$ if and only if
the corresponding solutions of the CYBE are equal.  
\end{thm}

\Pf . First the same argument as in the proof of Theorem
\ref{scthm} shows that $r(u,v)$ is uniquely determined
by terms $r_0$, $r_1$ and $r_2$. Furthermore, we have an equation
$$[r_2^{12}(v)-2r_2^{23}(v')-r_2^{13}(v+v')]\cdot u'-
[r_2^{12}(v)+r_2^{23}(v')+2r_2^{13}(v+v')]\cdot u=\ldots$$
where the RHS depends only on $r_0$ and $r_1$.
Hence, each term in the LHS can be recovered from $r_0$ and $r_1$.
Therefore, the same is true for the expression
$r_2^{12}(v)-r_2^{23}(v')$, hence for $r_2(v)$. 
The terms $r_0$ and $r_1$ are related by the equation
\begin{equation}\label{aux5}
r_0^{12}(v)r_0^{13}(v+v')-r_0^{23}(v')r_0^{12}(v)+r_0^{13}(v+v')r_0^{23}(v')=
r_1^{12}(v)+r_1^{23}(v')+r_1^{13}(v+v').
\end{equation}
We claim that $r_1$ is uniquely determined by $r_0$. Indeed,
let $v\mapsto s(v)$ be a $\Mat_n(\C)\otimes\Mat_n(\C)$-valued meromorphic
function in a neighborhood of zero such that $s^{21}(-v)=s(v)$ and
$$s^{12}(v)+s^{23}(v')+s^{13}(v+v')=0.$$
We have to prove that $s$ is zero. Applying $\pr\otimes\id\otimes\id$
to the equation we immediately deduce that
$$(\pr\otimes\id)(s(v))=0.$$
Similarly, $(\id\otimes\pr)(s(v))=0$, hence $s(v)=f(v)\cdot 1\otimes 1$
where $f(v)$ is an even meromorphic function satisfying 
$f(v)+f(v')+f(v+v')=0$. Hence, $f=0$.

It remains to show that $r_0(v)$ is uniquely determined by 
$\ov{r}(v)=(\pr\otimes\pr)(r_0(v))$ up to a summand of the form 
$cv\cdot 1\otimes 1$, where $c\in\C$, provided that $\ov{r}(v)$
has no infinitesimal symmetries. Let $(\wt{r}_0(v),\wt{r}_1(v))$
be another solution of the equation (\ref{aux5}) such that
$\wt{r}^{21}_0(-v)=-\wt{r}(v)$, $\wt{r}^{21}_1(-v)=\wt{r}_1(v)$.
We claim that if $(\pr\otimes\pr)(\wt{r}_0(v))=\ov{r}(v)$ then
$\wt{r}_0(v)=r_0(v)$. Indeed, we can write
$$\wt{r}_0(v)=r_0(v)+\phi^1(v)-\phi^2(-v)+\psi(v)\cdot 1\otimes 1$$
for some unique $\splin_n(\C)$-valued function $\phi(v)$ and some scalar
function $\psi(v)$. Let us denote the LHS of the equation (\ref{aux5})
by $LHS(r)$. Then we have
\begin{align*}
&0=(\pr\otimes\pr\otimes\pr)(LHS(\wt{r})-LHS(r))=\\
&\ov{r}^{12}(v)\cdot(\phi^3(-v')-\phi^3(-v-v'))+
\ov{r}^{23}(v')\cdot(\phi^1(v+v')-\phi^1(v))+
\ov{r}^{13}(v+v')\cdot(\phi^2(v')-\phi^2(-v)).
\end{align*}
If the function $\phi(v)$ is not constant then contracting this
equation with a generic functional in the third component we derive
that $\ov{r}(v)$ is a sum of two decomposable tensors which contradicts
our assumption. Hence, the function $\phi(v)$ has a constant value
$\phi\in\splin_n(\C)$. Now applying the projection
$\pr\otimes\pr\otimes\id$ to the difference of equations (\ref{aux5}) 
for $\wt{r}$ and $r$ we get the equation
\begin{equation}\label{aux6}
(\pr\otimes\pr\otimes\id)(\wt{r}_1^{12}(v)-r_1^{12}(v))=
(\pr\otimes\pr\otimes\id)(r_0^{12}(v)\phi^1-\phi^2r_0^{12}(v))-\phi^1\phi^2
+(\psi(v+v')-\psi(v'))\cdot\ov{r}^{12}(v)
\end{equation}
This is possible only if $\psi(v+v')-\psi(v')$ is independent of $v'$,
i.e. when $\psi$ is a linear function. Since $\psi(-v)=\psi(v)$ we
obtain $\psi(v)=cv$ for some constant $c\in\C$.
Thus, changing $r(u,v)$ to $\exp(cuv)r(u,v)$ we can assume that 
$\psi=0$. Finally making a substitution $v\mapsto -v$ and exchanging
the first two components in the equation (\ref{aux6}) we get (taking into
accound the unitarity condition) that
$$(\pr\otimes\pr)(r_0(v)\phi^1-\phi^2r_0(v))=
(\pr\otimes\pr)(-r_0(v)\phi^2+\phi^1r_0(v)),$$
or equivalently,
$$[\ov{r}(v),\phi^1+\phi^2]=0$$
which means that $\phi$ is an infinitesimal symmetry of $\ov{r}$.
Hence, $\phi=0$.
\ed

\begin{rems} 1. We don't know whether for every unitary non-degenerate
solution $\ov{r}(v)$ of the CYBE there exists a unitary solution of 
the AYBE of the form $\frac{1\otimes 1}{u}+r_0(v)+\ldots$ such that
$(\pr\otimes\pr)(r_0(v))=\ov{r}(v)$.

\noindent
2. In the case when $\ov{r}(v)$ has non-trivial
infinitesimal symmetries the proof
above shows that there are no more liftings of $\ov{r}(v)$ 
to a unitary solution
$r(u,v)$ of the AYBE (considered up to rescaling) than infinitesimal
symmetries of $\ov{r}(v)$. More precisely, such a lifting 
$r(u,v)=\frac{1\otimes 1}{u}+r_0(v)+\ldots$ is
uniquely determined by $r_0(v)$ and the difference between $r_0$'s for
two liftings always has form $\phi^1-\phi^2+c\cdot 1\otimes 1$
for some infinitesimal symmetry $\phi$ and some constant $c$.
\end{rems}

The above theorem can be applied in particular to the case when
$\ov{r}(v)$ is an elliptic non-degenerate solution of the CYBE.
Indeed, this follows from the following lemma (which I learned
from Pavel Etingof).

\begin{lem} Elliptic non-degenerate solutions of the CYBE
have no infinitesimal symmetries. 
\end{lem}

\Pf . The idea is to look at residues of such a solution at poles.
Let us denote $V=\C^n$.
Using the Killing form on $\splin(V)$ we can identify
$\splin(V)\otimes\splin(V)$ with endomorphisms of $\splin(V)$. 
Then the residues are operators corresponding to the action of the group
$G=(\Z/n\Z)^2$ on $\splin(V)$ induced by an irreducible projective
representation $\rho$ of $G$ on $V$ (see \cite{BD},5.1,5.2).
Let us denote by $\Ad\rho$ the representation of $G$ on $\splin(V)$.
It suffices to prove that if $A\in\SL(V)$ is such that
$$\Ad(A)\circ \Ad\rho(g)\circ\Ad(A)^{-1}=\Ad\rho(g)$$
then $A^{n^2}=1$. But this equation means that
for every $g\in G$ we have
$$A\rho(g)A^{-1}=c\cdot\rho(g)$$
for some constant $c\in\C^*$. Considering the determinants we see that
$c^n=1$, hence, 
$$A^n\rho(g)A^{-n}=\rho(g).$$
It follows that $A^n$ is scalar. Since it belongs to $\SL(V)$ we conclude
that $A^{n^2}=1$.
\ed

The conclusion one can draw from the above lemma and from theorem 
\ref{reconstr} is that elliptic solutions of the AYBE constructed
from triple Massey products on an elliptic curve can be uniquely
reconstructed from the limiting elliptic solutions of the CYBE.
As we have shown in \cite{P-hmc} the $A_{\infty}$-category
of elliptic curve (or at least the ``transversal'' part of it) can be
recovered from the usual category of vector bundles and from the
triple Massey products of the type considered in section \ref{Massey}. 
Hence, in some sense the information about all higher products
of the $A_{\infty}$-structure on elliptic curve (considered up to homotopy)
is encoded in elliptic solutions of the CYBE.

\bigskip

\section{Appendix}

In this appendix we prove two formulas for which we could not find 
references in the literature.
Let $\zeta(x,\tau)$ be the Weierstrass zeta-function associated
with the lattice $\Z+\Z\tau$. Let 
$\wp(x,\tau)=-\zeta'(x,\tau)$ be the corresponding $\wp$-function. 
For a pair of rational numbers $(r_1,r_2)$ we denote
\begin{equation}\label{zeta-char}
\zeta_{r_1,r_2}(x,\tau)=\zeta(x+r_1+r_2\tau,\tau)-r_1\eta_1(\tau)-
r_2\eta_2(\tau),
\end{equation}
where $\eta_1,\eta_2$ are quasi-periods corresponding to the basis $(1,\tau)$
(i.e. $\eta_1(\tau)=\zeta(x+1,\tau)-\zeta(x)$, $\eta_2(\tau)=
\zeta(x+\tau)-\zeta(x)$).
The first formula is
\begin{equation}\label{formula1}
\zeta(dx,d\tau)=\frac{1}{d}\cdot\sum_{i\in\Z/d\Z}\zeta_{\frac{i}{d},0}(x,\tau)+
\frac{x}{d}\cdot\sum_{i\in(\Z/d\Z)^*}\wp(\frac{i}{d},\tau).
\end{equation}
For the proof let us fix $\tau$ and
denote by $f(x)$ the difference between the LHS and the RHS.
Then one immediately checks that $f(x)$ is holomorphic on the entire plane,
$f'(x)$ is doubly periodic with respect to the lattice $\Z+\Z\tau$,
and $f(-x)=-f(x)$. Therefore, $f(x)=c\cdot x$ for some constant $c$.
Hence, it suffices to check the identity obtained from (\ref{formula1})
by differentiation:
$$\wp(dx,d\tau)=\frac{1}{d^2}\sum_{i\in\Z/d\Z}\wp(x+\frac{i}{d},\tau)-
\frac{1}{d^2}\sum_{i\in(\Z/d\Z)^*}\wp(\frac{i}{d},\tau).$$
But this can be proven directly from the definition of the $\wp$-function
as a series.

As a corollary of (\ref{formula1}) we immediately get that
$$\eta_2(d\tau)=\eta_2(\tau)+\frac{\tau}{d}\cdot\sum_{i\in(\Z/d\Z)^*}
\wp(\frac{i}{d},\tau).$$
Now it is easy to derive the following version of formula (\ref{formula1}):
\begin{equation}\label{formula1bis}
\zeta_{0,\frac{j}{d}}(dx,d\tau)=
\frac{1}{d}\cdot\sum_{i\in\Z/d\Z}\zeta_{\frac{i}{d},\frac{j}{d}}(x,\tau)+
\frac{x}{d}\cdot\sum_{i\in(\Z/d\Z)^*}\wp(\frac{i}{d},\tau).
\end{equation}

The second formula makes a connection between the special values of
the Kronecker function and Weierstrass zeta-function.
Namely using the notation (\ref{Kr-char}) we have 
\begin{equation}\label{formula2}
2\pi i F_{\frac{k}{d},\frac{l}{d}}(0,dx,d\tau)=
\sum_{j\in\Z/d\Z}\exp(-2\pi i\frac{kj}{d})
[\zeta_{\frac{j}{d},\frac{l}{d}}(x,\tau)-
\zeta_{\frac{j}{d},0}(-\frac{k\tau}{d},\tau)].
\end{equation}
where $d$, $k$ and $l$ are integers, $d>0$, $k$ is not divisible by $d$.
The proof of this formula is straightforward. Indeed, changing
$x$ one can reduce to the case $l=0$. Then the difference between the
LHS and the RHS is a holomorphic function of $x$, doubly periodic
with respect to the lattice $\Z+\Z\tau$, vanishing at $x=-\frac{k\tau}{d}$,
so it vanishes identically.

\end{document}